\documentclass{siamltex}

\usepackage{color}

\newtheorem{sideeg}           [theorem]{Example} 
\newenvironment{example}         {\begin{sideeg}\rm}{\end{sideeg}}
\newtheorem{sideremark}         [theorem]{Remark}
\newenvironment{remark}         {\begin{sideremark}\rm}{\end{sideremark}}

\newcommand{\BB}[1]{\mathbb{#1}}
\newcommand{\BS}[1]{\mathscr{#1}}

\newcommand{\ten}{\otimes}

\newcommand{\be}{\begin{equation}}
\newcommand{\ee}{\end{equation}}
\newcommand{\ba}{\begin{array}}
\newcommand{\ea}{\end{array}}
\newcommand{\la}{{\langle}}
\newcommand{\ra}{{\rangle}}
\newcommand{\R}                 {{\bf R}}
\newcommand{\C}                 {{\bf C}}

\newcommand{\N}                 {{\bf N}}

\newcommand{\er}[1]{\hbox{(\ref{#1})}}

\newcommand{\qed}               {\hfill $\square$}
\newcommand{\qH}{{\mathsf H}}
\newcommand{\qF}{{\mathsf F}}
\newcommand{\qh}{{\mathsf h}}
\newcommand{\cF}{{\mathcal F}}

\usepackage{amsmath}
\usepackage{amsfonts}
\usepackage{amssymb}
\usepackage{mathrsfs}
\usepackage{graphicx}

\title{An introduction to quantum filtering\thanks{
R.v.H.\ and L.B.\ are supported by the Army Research Office under Grant 
DAAD19-03-1-0073.  L.B.\ is additionally supported by the National 
Science Foundation under Grant PHY-0456720. M.R.J.\ is supported by the 
Australian Research Council.}}

\author{Luc Bouten\footnotemark[2] \and Ramon van Handel\footnotemark[2] 
\and Matthew R. James\footnotemark[3]}

\begin{document}
\maketitle

\renewcommand{\thefootnote}{\fnsymbol{footnote}}
\footnotetext[2]{
	Physical Measurement and Control 266-33, California Institute of
	Technology, Pasadena, CA 91125, USA (bouten@its.caltech.edu,
	ramon@its.caltech.edu).
}
\footnotetext[3]{
	Department of Engineering, Australian National University, 
	Canberra, ACT 0200, Australia (matthew.james@anu.edu.au).
}
\renewcommand{\thefootnote}{\arabic{footnote}}

\begin{abstract}
This paper provides an introduction to quantum filtering theory.
An introduction to quantum probability theory is given, focusing on the 
spectral theorem and the conditional expectation as a least squares 
estimate, and culminating in the construction of Wiener and Poisson 
processes on the Fock space.  We describe the quantum It\^o calculus and 
its use in the modelling of physical systems.  We use both reference 
probability and innovations methods to obtain quantum filtering equations 
for system-probe models from quantum optics.  
\end{abstract}



\pagestyle{myheadings}
\thispagestyle{plain}
\markboth{BOUTEN, VAN HANDEL, AND JAMES}{AN INTRODUCTION TO QUANTUM FILTERING}


\section{Introduction}
\label{sec intro}

Since even before the industrial revolution feedback control has played a
major role in the development of technology.  Nowadays many machines and
devices that make up our everyday lives use feedback to provide efficient
and reliable performance despite the ever increasing complexity and
miniaturization, and a rich control theory has been developed to aid in
the design of feedback controllers based on device models from classical
physics.  As microtechnology is making way for nanotechnology, however, we
are now rapidly approaching the boundary of the classical world past which
the effects of quantum mechanics cannot be neglected. 

The laws of quantum mechanics tell us that any description of the
phenomena at small scales is inherently nondeterministic in nature. This
opens new areas of application for stochastic control theory, which could
play an important role in a future generation of technology.  In
particular, as observations of quantum systems are inherently noisy, the
theory of filtering---the extraction of information from a noisy
signal---forms an essential part of any quantum feedback control strategy.

Quantum filtering was already implicit in early work on quantum
measurement theory by Davies in the 1960s \cite{Dav69,Dav76}. In its
modern form, the study of quantum filtering and control was pioneered by
Belavkin in a series of articles dating back to the early 1980s
\cite{Bel80,Bel83,Bel87,Bel92a,Bel92b}.  The theory developed by Belavkin
provides an essential foundation for statistical inference in e.g.\
quantum optical systems, and much of what we will discuss in the second
half of this article is based on his work.  The theory gained popularity
in the physics community after it was independently developed on a more
heuristic level by Carmichael in the early 1990s \cite{Car93} under the
name ``quantum trajectory theory'' and has since been widely applied in
the description of quantum optical experiments and as a computational
tool. 

Based on the foundations of quantum filtering theory, methods from
classical nonlinear and stochastic control can be developed and applied to
design feedback control laws for quantum systems. These methods may be
optimal in some sense, or otherwise designed with relevant considerations
in mind (e.g.\ stability). The resulting controllers are intended to be
implemented with some classical technology (e.g.\ digital or analog
electronics).  Recent experiments implementing quantum feedback controls
\cite{AASDM02,GSDM03,GSM04} have led to renewed interest in the field
which is now rapidly expanding \cite{Bel83,Wis94b,DJ99,DHJMT00,
BGM04,Jam04,Bou04b,HSM05,BEB05,Jam05,GBS05,EdB05,HaM05,HSM05a,BH05}.  We 
believe that a fruitful interaction between stochastic control and 
theoretical and experimental physics will be essential in paving the way 
towards the engineering of quantum technologies.

This paper provides an introduction to quantum filtering theory.  There
are three key ingredients that are required for the development of the
theory.  First, we need to capture both classical probability and quantum
mechanics within the framework of a generalized probability theory, called
noncommutative or quantum probability theory.  The central object in this
theory, the spectral theorem, provides a link between quantum systems and
the associated probabilistic measurement outcomes.  Second, we need a
noncommutative generalization of the concept of conditional expectations.  
As in classical probability, we will find that a suitably restricted
definition of the quantum conditional expectation is none other than a
least squares estimator, which elucidates its role in quantum filtering
theory.  Finally, we need a noncommutative analog of stochastic calculus
and quantum stochastic differential equations (QSDE).  This provides a
broad class of models for which we can obtain filtering equations.

A typical physical scenario, to which the theory that we will develop can
be applied, is illustrated schematically in Figure \ref{fig:model}.  A
cloud of (usually cold, trapped) atoms interacts with the electromagnetic
field in free space; this can be coherent light from a laser, or even the
vacuum.  Depending on their internal state the atoms can, for example,
emit radiation into the field.  If we detect this radiation using an
optical detection setup we can try to infer some information on the
internal state of the atoms---this is precisely the goal of quantum
filtering theory.  If we wanted to control the state of the atoms, we
could then feed back some function of the state estimates through a
suitable actuator.  Recent laboratory experiments, e.g.\ \cite{GSM04},
implement precisely such a setup, and provide a motivating example for the
theory.

We begin in \S\ref{sec BG} by providing some background for quantum
filtering. This includes a discussion of the quantum mechanics and quantum
probability in the simplest, finite-dimensional context. In \S\ref{sec QP}
quantum probability is developed in detail. Then in \S\ref{sec Fock} we
show how Wiener and Poisson processes emerge in a particular quantum
probabilistic model based on the Fock space, and how these can be used to
develop a noncommutative stochastic calculus. In \S\ref{sec filtering} we
introduce a class of system-observation models that describe typical
experiments in quantum optics. \S\ref{sec KS formula} deals with the
derivation of quantum filtering equations using the reference probability
approach, while \S\ref{sec mtg} gives an alternative derivation using the
innovations or martingale method.

{\bf Scope.} It has been our aim to make quantum probability and filtering
theory accessible, modulo a set of technicalities, to readers with a
minimal number of prerequisites.  We (only) presume some familiarity with
probability theory and elementary functional analysis.  We have put an
emphasis on introducing the mathematical structures of quantum probability
theory and on demonstrating their significance and their use.  As a
consequence we do not everywhere achieve the highest level of rigor; we
are particularly lax in the use of unbounded operators and their domains.
It is our hope that skimming over these technicalities has enabled us to
paint a clearer picture of the pillars of the theory and of the essential
techniques involved.  That being said, we should point out that many of
the tools described in this paper are applied regularly and successfully
by physicists without paying any attention to the technical issues
involved; the reader should not hesitate to get his feet wet!

It is an ambitious project to introduce an unfamiliar probability theory, 
a new stochastic calculus, and to even solve a nontrivial problem 
(filtering) within the confines of about 40 pages.  Though we have tried 
to give a pedagogical treatment, the explanations are sometimes 
necessarily terse; we hope that the reader will be sufficiently compelled 
to work his way through the paper.  Needless to say there are many 
omissions; one that particularly deserves mention is the linear case: 
indeed, the quantum Kalman filter, and the corresponding theory of quantum 
LQG control, can be developed along similar lines to the filters we will 
discuss.  We have chosen to omit this topic in order to avoid the 
technicalities of QSDEs with unbounded coefficients, but refer instead 
to \cite{EdB05} and the references therein.

{\bf Notation.} The sets of natural, real and complex numbers are denoted
$\N$, $\R$ and $\C$ respectively. In general, script symbols (e.g.\
$\mathscr{Y}$) are used for von Neumann algebras, while calligraphic
symbols (e.g.\ $\mathcal{Y}$) stand for $\sigma$-algebras.  $\mathcal{B}$
is the Borel $\sigma$-algebra on $\R$.  Classical probability spaces are
denoted as $(\Omega,\mathcal{F},{\bf P})$, and $E_{\bf P}$ denotes the
expectation with respect to the measure ${\bf P}$. Blackboard symbols
(e.g.\ $\mathbb{P}$) denote states on von Neumann algebras.  Sans-serif
symbols (e.g.\ $\qH$) are used for Hilbert spaces.  Hilbert space
adjoints, as well as the scalar complex conjugate, are indicated by $^*$,
and the Hilbert space inner product is denoted by
$\langle\cdot,\cdot\rangle$.  The commutator of two bounded operators is
denoted by $[X,Y]=XY-YX$.  $I$ is the identity operator.

 
\section{Background and motivation}
\label{sec BG}

In this article we adopt a modern quantum probability formulation of
quantum mechanics.  Quantum probability is the noncommutative counterpart
of Kolmogorov's axiomatic characterization of classical probability
theory.  In addition to the natural interpretation and mathematical tools
provided by Kolmogorov's formalism, one of its major successes is that
conditioning is a derived concept rather than an additional axiom.  The
situation is much the same in quantum probability; in particular, the
conditioning axiom or ``projection postulate'' as it is traditionally
posed in quantum mechanics can emerge as a consequence of conditional
expectation and the physical idea that in a single experiment one only has
direct access to information contained in a commutative subalgebra of
observables.

Considering the success of the classical (Kolmogorov) theory, it should
come as no surprise that the mathematical abstraction provided by the
framework of quantum probability pays off significantly (as we will see
throughout the article).  Introductory physics textbooks on quantum
mechanics rarely use such a description, however. In this section we
introduce the basic concepts of quantum probability in their simplest 
form, and attempt to provide contact with ideas about quantum mechanics 
that readers may be familiar with. This is intended to provide a reference 
point for interpreting the quantum probabilistic framework used in this 
paper.

\subsection{Some textbook quantum mechanics}
\label{sec BG QM}

According to the textbook by Merzbacher \cite[page 1]{EM98}, ``Quantum
mechanics is the theoretical framework within which it has been found
possible to describe, correlate, and predict the behavior of a vast range
of physical systems, from particles through nuclei, atoms and radiation to
molecules and condensed matter.'' Central to quantum mechanics are the
notions of {\em observables}, which are mathematical representations of
physical quantities that can (in principle) be measured, and {\em states},
which summarize the status of physical systems and permit the calculation
of statistical quantities (such as probabilities, expectations,
correlations) of observables. 

Indeed, the reader may be familiar with the {\em Schr{\"o}dinger
wavefunction} $\psi(q,t)$ for a particle of mass $m$ moving in a force
field $V(q)$ (dependent on position $q$, in one dimension for simplicity). 
If $Q$ is the observable representing position (defined below in Example
\ref{eg q and p}), the expected position of the particle when in a state
described by $\psi(q,t)$ at time $t$ is defined to be
\be
	 \la Q \ra = \int  q \vert \psi(q,t) \vert^2 dq .
	 \label{bg-expect-1}
\ee
The wavefunctions are normalized to one $ \int  \vert \psi(q,t) \vert^2 dq =1$,
so that $\vert \psi(q,t) \vert^2$ could be interpreted as the probability
density of the position of the particle.  The dynamics of the particle are
described by the famous {\em Schr{\"o}dinger equation} 
\be
 	 i \hbar \frac{\partial \psi(q,t)}{\partial t} = - \frac{\hbar^2}{2m} 
	 \frac{\partial^2 \psi(q,t)}{\partial q^2} + V(q) \psi(q,t) ,
	 \label{bg-schrodinger}
\ee
where $\hbar=h/2\pi$, $h$ is {\em Planck's constant}, and $i^2=-1$.
 
The key distinction between classical (i.e.\ non-quantum) and quantum 
mechanics is that quantum mechanics is {\em noncommutative}, meaning that
there exist observables that do not commute, a fact which has deep 
implications.  The momentum observable $P$ (defined below in Example 
\ref{eg q and p}) does not commute with the position observable $Q$; in 
fact $[Q, P] = QP-PQ = i \hbar \,I$.  The most famous implication of this 
failure of commutativity is {\em Heisenberg's uncertainty relation}, 
which asserts that 
\be
	\Delta Q \Delta P \geq \frac{1}{2} \vert \la i [Q,P ] \ra \vert = 
	\frac{\hbar}{2}, \label{hu-1}
\ee
where the variances are defined by
$\Delta Q = (\la Q^2 \ra - \la Q \ra^2)^{1/2}$, 
$\Delta P =(\la P^2 \ra - \la P \ra^2)^{1/2}$.
Naive interpretation of the Heisenberg uncertainty relation can be 
misleading; we will discuss its precise meaning in the following section.
Nonetheless, it evidently implies that there is a fundamental irreducible 
randomness in quantum mechanics. This is in contrast to classical 
randomness, which in principle can be eliminated with enough effort and
information. Experimental evidence has repeatedly confirmed the 
irreducible randomness of quantum mechanical observations.

Let us make this somewhat vague discussion a little more precise.  For 
simplicity, we will work in this section only in a finite-dimensional 
setting (in which observations can only take a finite number of values, 
i.e.\ they are finite-state random variables).  First, recall that if 
$A=A^*$ is a self-adjoint operator on a finite dimensional Hilbert space 
$\qH=\C^n$, it has at most $n$ (distinct) real eigenvalues.  The set  
${\rm spec}(A)=\{ a_j \}$ of eigenvalues of $A$ is called the spectrum 
of $A$, and $A$ can be written as 
\be
	A = \sum_{a\in{\rm spec}(A)} a\,P_a ,
	\label{spectral-1}
\ee
where $P_a$ is the projection operator onto the subspace of $\qH$ spanned by
vectors with eigenvalue $a$.  The projections resolve the identity $\sum_{a
\in {\rm spec}(A)} P_a = I$.

In this finite-dimensional setting, the following operational
characterization of quantum mechanical models (often referred to as the 
``postulates'' of quantum mechanics) can be found in most introductory 
textbooks.

{\bf Observables}.
Physical quantities like position, momentum, spin, etc., are represented by
self-adjoint operators on the Hilbert space $\qH$ and are called {\em
observables}.   These are the noncommutative counterparts of random 
variables.

{\bf States}.
A state is meant to provide a summary of the status of a physical system
that enables the calculation of statistical quantities associated with
observables.  A generic state is specified by a {\em density matrix} 
$\rho$, which is a self-adjoint operator on $\qH$ that is positive 
$\rho\ge 0$ and normalized ${\rm Tr}[\rho]=1$.  This is the noncommutative 
counterpart of a probability density.

{\bf Measurement}.
A {\em measurement} is a physical procedure or experiment that produces
numerical results related to observables. In any given measurement, the
allowable results take values in the spectrum ${\rm spec}(A)$ of a chosen
observable $A$.   Given the state $\rho$, the value $a\in{\rm spec}(A)$ 
is observed with probability ${\rm Tr}[\rho P_a]$.  Consequently, the 
expectation of an observable $A$ is given by $\langle A\rangle={\rm 
Tr}[\rho A]$.

{\bf Conditioning}.
Suppose that a measurement of $A$ gives rise to the observation $a\in{\rm 
spec}(A)$.  Then we must condition the state in order to predict the outcomes
of subsequent measurements, by updating the density matrix $\rho$ using
\be
	\rho\mapsto\rho'[a]=\frac{P_a\rho P_a}{{\rm Tr}[\rho P_a]}.
	\label{axiom-proj}
\ee
This is known as the ``projection postulate''.

{\bf Evolution}.
A {\em closed} (i.e.\ isolated) quantum system evolves in a {\em unitary} 
fashion: a physical quantity that is described at time $t=0$ by an 
observable $A$ is described at time $t>0$ by $A(t)=U(t)^*AU(t)$, where 
$U(t)$ is a unitary operator for each time $t$.  The unitary is generated 
by the {\em Schr{\"o}dinger equation}
\be
	i \hbar \frac{d}{dt}U(t) =  H(t)U(t) ,
	\label{schrodinger}
\ee
where the (time dependent) Hamiltonian $H(t)$ is a self-adjoint operator 
for each $t$.

Before continuing, we make the following remarks. 
\begin{remark}(Pure states).
The set of density matrices $\rho$ is convex; we can thus wonder what are
the extremal points in this set, i.e.\ those that correspond to the most
informative states.  It is not difficult to show that the set of extremal
density matrices is the set of projections onto one-dimensional subspaces.
Thus we can specify any extremal state uniquely (up to a phase factor
$e^{i\varphi}$) by a single unit vector $\psi\in\qH$ in the corresponding
subspace, and ${\rm Tr}[\rho X]=\langle\psi,X\psi\rangle$ for any operator
$X$.  In classical probability theory, the set of probability measures is 
also convex and the extremal measures are deterministic (Dirac) measures.  In
the quantum mechanical setting, on the other hand, the Heisenberg
uncertainty relation implies that even extremal states do not give
deterministic measurement outcomes for all observables. 

Historically, and in most textbooks,  quantum mechanics is first 
formulated in terms of the extremal states (called {\em pure states}) and 
the description is later generalized to density matrices ({\em mixed 
states}).  The Schr{\"o}dinger wavefunction $\psi(q,t)$ is an example of a 
pure state vector in an infinite-dimensional Hilbert space setting. \qed
\end{remark}
\begin{remark}(Heisenberg vs.\ Schr{\"o}dinger picture).
In the above description of time evolution we work with a fixed state 
while the observables change in time.  This conforms to the usual 
treatment in classical probability theory, where the underlying 
probability measure is fixed at the outset and the random variables are 
time dependent (stochastic processes).  In quantum mechanics this is 
known as the {\em Heisenberg picture}; equally (or perhaps more) popular
is the {\em Schr{\"o}dinger picture}, in which the observables are 
considered fixed and the density matrix evolves as $\rho(t)=U(t)\rho U(t)^*$.
The two pictures are essentially equivalent as ${\rm Tr}[\rho A(t)]={\rm 
Tr}[\rho(t)A]$ for any observable $A$.

Note that if we start in a pure state, then unitary evolution preserves 
this property; in terms of the state vector, $\psi(t)=U(t)\psi$.
Intuitively, this enforces the physical idea that no information is lost
from an isolated system.  Together with (\ref{schrodinger}) we obtain the 
traditional Schr{\"o}dinger equation for $\psi(t)$, of which
(\ref{bg-schrodinger}) is a special case (for a specific choice of $H$, in
infinite dimensions).  We will always work in the Heisenberg picture,
however, as we will be dealing with (quantum) stochastic processes. \qed
\end{remark}

As a basic illustration we discuss the following simple example.

\begin{example} \label{eg SG 1} 
One of the classic experimental demonstrations of the necessity of quantum
mechanics was performed in 1922 by Stern and Gerlach.  A silver atom is
subjected to an inhomogeneous magnetic field.  The atom possesses an
intrinsic magnetic moment, and hence experiences a force that is
proportional to the component of its magnetic moment in the direction of
the field gradient.  As Stern and Gelach did not prepare the atom in a
particular orientation, they expected it to be deflected randomly in a 
continuous range of directions corresponding to a random orientation of the 
magnetic moment.  Repeated runs of the experiment showed, however, that 
the atom is randomly deflected into two discrete directions only---the 
reason being that in quantum mechanics the magnetic moment (or spin) 
observable is discrete, rather than continuous.  Atoms deflected in the 
upper direction are said to have \lq\lq{spin up}\rq\rq, while those in 
the lower direction have \lq\lq{spin down}\rq\rq. 

A simple model of a spin is as follows.  Let $\qH=\C^2$, and consider the 
observable
\be
	\sigma_z = \left( \ba{cc} 1 & 0 \\ 0 & -1 \ea \right)
	\label{bg-sigma-z}
\ee
representing spin in the $z$ direction.   We have 
${\rm spec}(\sigma_z)=\{-1,1\}$, which correspond to spin down and spin 
up, respectively.  In terms of the eigenprojections
$$
	P_{z,1} = \left( \ba{cc} 1 & 0 \\ 0 & 0 \ea \right), \qquad
	P_{z,-1} = \left( \ba{cc} 0 & 0 \\ 0 & 1 \ea \right),
$$
we can write $\sigma_z=P_{z,1}-P_{z,-1}$.  The next step is to introduce a 
state.  Consider a pure state, given by the vector $\psi=(c_1~c_{-1})^T$ 
with $|c_{1}|^2+|c_{-1}|^2=1$.  If we observe $\sigma_z$, we obtain the 
outcome $1$ (spin up) with probability $\langle\psi,P_{z,1}\psi\rangle=
|c_1|^2$, or the outcome $-1$ with probability 
$\langle\psi,P_{z,-1}\psi\rangle=|c_{-1}|^2$.  \qed
\end{example}

\subsection{A first look at quantum probability}
\label{sec BG towards}

The description of quantum mechanics in the previous section contains the
rudiments of a viable probability theory.  We will now formalize these
ideas, once again restricting ourselves to the finite-dimensional case for
simplicity (the general theory, which will be discussed in \S\ref{sec QP},
is conceptually very similar). Two key ideas, which we elaborate on below,
form the essence of the formalism: the first is that a set of measurements
made in a single realization\footnote{
	By a {\it realization} or an {\it experiment} we mean that
	random variables are assigned a definite value, as is the
	case if we perform measurements on a single physical system.
	In classical probability this corresponds to the choice
	of a sample point $\omega\in\Omega$; the quantum case
	is a little more subtle.
} of a quantum experiment corresponds to a particular choice of a
commutative algebra of observables; and the second is that any such
commutative algebra is entirely equivalent to a classical (Kolmogorov)
probability model.

A classical probability model is described by a probability space
$(\Omega, {\mathcal F}, {\mathbf P})$.  Here $\Omega$, the sample space,
is not of essential importance; the basic ingredients of the theory are
the events that can occur, contained in the $\sigma$-algebra
$\mathcal{F}$, and their probabilities, which are determined by the
measure ${\bf P}$.  Equivalently, we could describe an event
$F\in\mathcal{F}$ by a random variable $\chi_F$ which takes the value $1$
if $F$ occurs and $0$ otherwise (the indicator function on $F$), and the
probability of the event is simply the expectation of $\chi_F$.  We have
already encountered such objects in the previous section: events are
precisely those observables that are projection operators ($P=P^*=P^2$),
and the probability of an event $P$ is given by $\mathbb{P}(P)={\rm
Tr}[\rho P]$.  Thus the set of projections, together with the linear map
$\mathbb{P}$, play much the same role as the classical pair
$\mathcal{F},{\bf P}$. 

We run into trouble in the quantum case when we try to ascribe joint
probabilities to certain events.  This is always possible in classical
probability theory: the joint probability of the events $A$ and $B$ is
${\bf P}(A\cap B)=E_{\bf P}(\chi_A\chi_B)$.  But given two projection
operators $P,Q$ the operator $PQ$ is not guaranteed to be a projection or
even an observable ($(PQ)^*=QP$), {\em unless} $P$ and $Q$ commute.  This
simple observation is no coincidence; it has the following physical
interpretation: in a single realization of a quantum probability model, we
can only verify the truth of a set of commuting events.  This is in
contrast with classical probability where in every realization any event
is either true or false, whether we choose to observe it or not.  In
quantum probability we can a priori choose to verify the truth of an
arbitrary event, but subsequently some of the other events (those that do
not commute with the observed event, said to be {\em incompatible}) become 
meaningless within the same realization. 

The incompatibility of events is a significant conceptual departure from 
classical probability, and requires a little getting used to.  In many 
ways, however, this is the only essential departure from classical 
probability theory.  We now begin to contruct the mathematical formalism 
of quantum probability, and we will show that it is indeed very close to 
Kolmogorov's theory.

Consider the following idea.  Suppose we decide to measure an observable
$A$ and obtain a particular outcome $a\in{\rm spec}(A)$.  Then we do not
need to perform another measurement to know that any function $f(A)$ would
give the outcome $f(a)$; in essence, this is merely a relabeling of the
measurement outcomes of $A$.  Indeed,
\be
	A=\sum_{a\in{\rm spec}(A)}a\,P_a\quad
	\Longrightarrow\quad
	f(A)=\sum_{a\in{\rm spec}(A)}f(a)\,P_a,
\ee
and all such operators commute with each other. Thus measuring $A$
``automatically'' measures all functions $f(A)$. The set of operators
$\mathscr{A}=\{X:X=f(A),~f:{\bf R}\to{\bf C}\}$ forms a commutative
$*$-algebra, i.e.\ arbitrary (complex) linear combinations, products and
adjoints of operators in $\mathscr{A}$ are still in $\mathscr{A}$,
$I\in\mathscr{A}$, and all elements of $\mathscr{A}$ commute. We will call
$\mathscr{A}$ the $*$-algebra generated by $A$\footnote{
	In fact, it is the smallest $*$-algebra of operators that contains $A$.
}.  A linear map
$\mathbb{P}:\mathscr{A}\to{\bf C}$ that is positive ($\mathbb{P}(A)\ge 0$
if $A\ge 0$) and normalized ($\mathbb{P}(I)=1$) is called a state on
$\mathscr{A}$ (clearly we can always write such a state as $A\mapsto{\rm
Tr}[\rho A]$ for some density matrix $\rho$). Note that the projections
$P\in\mathscr{A}$ are precisely those events that we can distinguish by
measuring $A$, and $\mathbb{P}(P)$ gives their probabilities.   
We can similarly generate the commutative $*$-algebra of functions 
of an arbitrary set of commuting observables.

The algebraic structure we have introduced is of fundamental importance
as it provides us with a direct connection to the classical theory, as 
follows:

\begin{theorem}[Spectral theorem, finite-dimensional case]
\label{thm spectral finite}
Let $\mathscr{A}$ be a commutative $*$-algebra of operators on a 
finite-dimensional Hilbert space, and let $\mathbb{P}$ be a state on 
$\mathscr{A}$. Then there is a probability space $(\Omega,\mathcal{F},{\bf P})$
and a map $\iota$ from $\mathscr{A}$ onto the set of measurable functions 
on $\Omega$ that is a $*$-isomorphism, i.e.\ a linear bijection with 
$\iota(AB)=\iota(A)\iota(B)$ (pointwise) and $\iota(A^*)=\iota(A)^*$, and 
moreover $\mathbb{P}(A)=E_{\bf P}(\iota(A))$.
\end{theorem}

\begin{proof} 
The proof is an elementary exercise in linear algebra.  As the Hilbert
space $\mathsf{H}$ has dimension $n<\infty$, we can without loss of
generality suppose that $\mathsf{H}={\bf C}^n$ and that $\mathscr{A}$ is a
commutative $*$-algebra of complex $n\times n$ matrices.  As all the
elements of $\mathscr{A}$ commute, we can find a unitary matrix $U$ such
that $U^*AU$ is a diagonal matrix for every $A\in\mathscr{A}$.  Let
$\Omega=\{1,\ldots,n\}$.  Define $\iota(A):\Omega\to{\bf C}$ by
$\iota(A)(i)=(U^*AU)_{ii}$ for every $A\in\mathscr{A}$.  Next, define
$\mathcal{F}=\sigma\{\iota(A):A\in\mathscr{A}\}$.  Finally, define ${\bf
P}(S)=\mathbb{P}(\iota^{-1}(\chi_S))$ for every $S\in\mathcal{F}$.  We 
have now explicitly constructed $(\Omega,\mathcal{F},{\bf P})$ and $\iota$.
\qquad 
\end{proof}

Evidently the commutative $*$-algebra structure is completely equivalent to 
classical probability theory; by simultaneously diagonalizing all the 
operators in the algebra, we obtain an explicit representation of 
measurable random variables as the functions on the diagonals.  We also 
note the following.  Suppose we are given some (large) commutative 
$*$-algebra $\mathscr{A}$, and consider a subalgebra 
$\mathscr{B}\subset\mathscr{A}$ generated by a single element 
$B\in\mathscr{A}$.  If we apply the map $\iota$ to $\mathscr{B}$, we obtain 
precisely the subset of functions on $\Omega$ that are measurable with 
respect to $\sigma\{\iota(B)\}$.  Thus subalgebras play the same role in 
quantum probability as sub-$\sigma$-algebras in classical probability; 
they allow us to keep track of particular subsets of information.

We do not a priori have a basis for specifying a particular commutative 
$*$-algebra; given a quantum system, we could decide to measure any of a 
large set of incompatible observables.  The discussion up to this point 
motivates the following definition.

\begin{definition}[Quantum probability space, finite-dimensional case] 
\label{def fd-qps}
	A pair $({\mathscr N},\mathbb{P})$, where ${\mathscr N}$ is a 
	$\ast$-algebra of operators on a finite-dimensional Hilbert space and 
	$\mathbb{P}$ is a state on $\mathscr{N}$, is called a 
	(finite-dimensional) quantum probability space.
\end{definition}

Usually we will choose $\mathscr{N}$ to be the set of all (bounded)
operators $\mathscr{B}(\mathsf{H})$ on some underlying Hilbert space
$\mathsf{H}$.  The principles of quantum probability now boil down to the
following.  In each realization, we must make a choice of commutative
$*$-subalgebra $\mathscr{A}\subset\mathscr{N}$ which fixes the
observations.  Every statistic that pertains to these observations (e.g.,
the statistics compiled by repeating the experiment many times with the
same choice of $\mathscr{A}$) is now described by the classical
probability model obtained through the spectral theorem.  The reader
should convince himself that the operational description given in the
previous section fits neatly within this model (with the exception of
conditioning, which we discuss \S\ref{sec BG condex}).

Notice that in contrast to a classical probability space $(\Omega, 
{\mathcal F}, {\mathbf P})$, there are no sample points $\omega \in \Omega$  
in a quantum probability space.  The sample points emerge through the 
spectral theorem after the choice of a commutative $*$-subalgebra.

\begin{example}
\label{eg SG 2}
Let us reformulate Example \ref{eg SG 1}.  Set $\qH = \C^2$ and choose
${\mathscr N}={\mathscr B}(\qH)=M_2$, the $\ast$-algebra of $2\times 2$
complex matrices.  The pure state is defined by
$\mathbb{P}(A)=\langle\psi,A\psi\rangle=\psi^*A\psi$ (recall that
$\psi=(c_1~c_{-1})^T$ with $|c_1|^2+|c_{-1}|^2=1$). 

The observable $\sigma_z$, used to represent spin measurement in the $z$
direction, generates a commutative $*$-subalgebra ${\mathscr A}_z \subset
{\mathscr N}$.  It is not difficult to see that ${\mathscr A}_z$ is simply
the linear span of the events $P_{z,1}$ and $P_{z,-1}$. Let us now apply
the spectral theorem; we obtain the probability space
$(\Omega,\mathcal{F},{\bf P})$ where $\Omega=\{1,2\}$,
$\mathcal{F}=\{\varnothing,\{1\},\{2\},\Omega\}$, ${\bf
P}(\{1\})=|c_1|^2$, etc., and $\iota(P_{z,1})=\chi_{\{1\}}$,
$\iota(P_{z,-1})=\chi_{\{2\}}$.  In particular, the random variable
$\iota(\sigma_z):(1,2)\mapsto(1,-1)$ has precisely the right properties. 

Now suppose we do not wish to measure the intrinsic angular momentum (spin)
in the $z$-direction, but in the $x$-direction.   This corresponds to the 
observable
\be
	\sigma_x = \left( \ba{cc} 0 & 1 \\ 1 & 0 \ea \right),
	\label{bg-sigma-x}
\ee
which has the spectral decomposition $\sigma_x=P_{x,1}-P_{x,-1}$ with
$$
	P_{x,1} = \frac{1}{2}
	\left( \ba{cc} 1 & 1 \\ 1 & 1 \ea \right), \qquad
	P_{x,-1} = \frac{1}{2}
	\left( \ba{cc} 1 & -1 \\ -1 & 1 \ea \right).
$$
The observable $\sigma_x$ also generates a commutative $*$-subalgebra 
$\mathscr{A}_x={\rm span}\{P_{x,1},P_{x,-1}\}$ to which we can apply the 
spectral theorem.  However, as $\sigma_x$ and $\sigma_z$ do not commute, 
they cannot be jointly represented on a classical probability space 
through the spectral theorem.  In other words, $\sigma_x$ and $\sigma_z$ are 
incompatible and their joint statistics are undefined; hence they cannot 
both be observed in the same realization.
\qed
\end{example}
 
To conclude this section, let us say a few words about the interpretation
of the Heisenberg uncertainty relation.  The relation says that the
product of the variances of two noncommuting observables is bounded from
below by a positive constant.  It is important to realize, however, that
the two observables cannot be measured in the same realization as they are
incompatible---in particular, the {\em co}variance of the observables is
undefined.  Rather, the uncertainty relation is a statement about the
properties of quantum states: for any state, the statistics of the two 
observables, compiled in the course of separate realizations in each of 
which only one of the observables is measured, must obey the Heisenberg 
inequality\footnote{
	In the physics literature one often find statements to the effect
	that the Heisenberg uncertainty relation limits the precision with which 
	we can ``imperfectly'' observe two noncommuting observables 
	simultaneously, i.e.\ within the same realization.  This is a 
	misconception.  Though the idea of an imperfect measurement can be 
	implemented rigorously (e.g.\ \cite{Hol82}), this gives rise to an 
	uncertainty relation which is different than Heisenberg's
	uncertainty relation \cite{ArthursKelly}. 
}.

\subsection{Composite systems}
\label{sec BG composite}

We will often wish to form a composite probability model from two 
separate probability spaces.  In classical probability theory, two 
probability spaces $(\Omega_1,\mathcal{F}_1,{\bf P}_1)$ and
$(\Omega_2,\mathcal{F}_2,{\bf P}_2)$ can be merged into a single 
probability space 
$(\Omega_1\times\Omega_2,\mathcal{F}_1\times\mathcal{F}_2,{\bf 
P}_1\times{\bf P}_2)$ where ${\bf P}_1\times{\bf P}_2$ is the product 
measure.  We now briefly describe the noncommutative counterpart.

Consider a composite system constructed from two quantum probability spaces
$(\mathscr{N}_1,\mathbb{P}_1)$, $(\mathscr{N}_2,\mathbb{P}_2)$ of 
operators on the Hilbert spaces $\qH_1$ and $\qH_2$, respectively.
The composite quantum probability space consists of operators on the 
tensor product Hilbert space $\qH_1\otimes\qH_2$; for vectors 
$\psi_1,\phi_1\in\qH_1$ and $\psi_2,\phi_2\in\qH_2$, the inner product on 
$\qH_1\otimes\qH_2$ is given by
$$
	\la \psi_1 \otimes \psi_2, \phi_1 \otimes \phi_2 \ra = \la \psi_1,
	\phi_1 \ra \la \psi_2, \phi_2 \ra,
$$
which is extended by linearity to any vector in the tensor product space.
The algebra $\mathscr{N}_1\otimes\mathscr{N}_2$ is generated by elements of 
the form
$$
	(A_1\otimes A_2)(\psi_1\otimes\psi_2)=A_1\psi_1\otimes A_2\psi_2,
$$
where $A_1\in\mathscr{N}_1$ and $A_2\in\mathscr{N}_2$.  Finally, the 
product state is defined by
$$
	(\mathbb{P}_1 \otimes \mathbb{P}_2) (A_1 \otimes A_2 )  
	= \mathbb{P}_1(A_1)   \mathbb{P}_2 (A_2 ),
$$
and is extended by linearity.  The quantum probability space 
$(\mathscr{N}_1\otimes\mathscr{N}_2,\mathbb{P}_1 \otimes \mathbb{P}_2)$
of operators on the Hilbert space $\qH_1\otimes\qH_2$ describes the 
composite system.  The reader should verify that if $\mathscr{N}_1$ and 
$\mathscr{N}_2$ are commutative, then applying the spectral theorem to 
the composite system is equivalent to applying the spectral theorem to 
the individual subsystems, then forming the composite classical 
probability space.

\subsection{Conditional expectations}
\label{sec BG condex}

Let us recall for a moment the Stern-Gerlach experiment of Examples
\ref{eg SG 1} and \ref{eg SG 2}.  We have introduced the observables
$\sigma_z$ and $\sigma_x$, corresponding to spin in the $z$ and $x$
directions.  These observables are incompatible, so we cannot measure them
in the same realization.  Recall that in order to measure $\sigma_z$,
Stern and Gerlach apply a field gradient in the $z$ direction; the atom
then acquires momentum in that direction proportional to $\sigma_z$, and
we can determine the value of $\sigma_z$ in that realization by observing
whether the atom is deflected up ($1$) or down ($-1$).  Similarly,
$\sigma_x$ is measured by orienting the field gradient along the $x$ axis. 

We wouldn't be measuring both $\sigma_z$ and $\sigma_x$ by applying both
field gradients simultaneously: rather, as magnetic fields add
vectorially, this would measure the spin in some other direction in the
$x$-$z$ plane whose observable commutes with neither $\sigma_z$ nor
$\sigma_x$.  On the other hand, we could first apply the field gradient in
the $z$ direction until we can resolve $\sigma_z$, then turn this field
off and switch on a field in the $x$ direction to resolve $\sigma_x$.  It
is a characteristic feature of quantum mechanics that the measurement
outcomes in such a procedure can differ drastically depending on what
order we apply the fields.  It is thus of crucial importance to specify
precisely how such measurements are performed by including in the quantum
probability space a model of the measurement apparatus (or {\em probe}). 

We defer the discussion of the Stern-Gerlach measurement with magnetic
fields until we have developed the necessary machinery in \S\ref{sec QP}. 
For sake of example, we develop in this section a simpler probe model
which shows the main features of the procedure.  We will see that this
probe model, together with the concept of conditional expectations,
reproduces precisely the traditional projection postulate of \S\ref{sec BG
QM}. 

Let us begin by discussing conditional expectations in the noncommutative
context.  The key observation we need is the following.  The conditional
probability of an event $B$ given an event $A$ is the probability that $B$
is true given that $A$ is true in the same realization.  Hence the concept
of conditioning inherently makes sense only in the context of quantities
that can be observed in the same realization of an experiment. This means
that we can only define conditional expectations in commutative
subalgebras of a quantum probability space; but as long as we are 
restricted to the commutative case, the spectral theorem allows us to 
define any probabilistic operation directly in terms of the associated 
classical probability space (see \cite{BH05}).

To be more precise, let $(\mathscr{N},\mathbb{P})$ be a quantum
probability space, $\mathscr{A}\subset\mathscr{N}$ a commutative
subalgebra and $B\in\mathscr{N}$ a self-adjoint element commutes with
every $A\in\mathscr{A}$.  Then $B$ and $\mathscr{A}$ generate a larger
commutative subalgebra $\mathscr{C}\subset\mathscr{N}$, to which we can
apply the spectral theorem to obtain a $*$-isomorphism $\iota$.  The
conditional expectation is now simply inherited from the classical space
as $\mathbb{P}(B|\mathscr{A})= \iota^{-1}(E_{\bf
P}(\iota(B)|\sigma\{\iota(\mathscr{A})\}))$.  Note, however, that if $B,C$ 
are two self-adjoint operators that commute with every $A\in\mathscr{A}$, 
this does not necessarily imply that $B$ and $C$ commute.  The set
$\mathscr{A}'=\{B\in\mathscr{N}:AB=BA~\forall A\in\mathscr{A}\}$,
the {\em commutant} of $\mathscr{A}$ (in $\mathscr{N}$), is the largest
$*$-subalgebra of operators that can be conditioned on $\mathscr{A}$.  
The conditional expectation is defined as above for its self-adjoint 
elements, and extends to all of $\mathscr{A}'$ by linearity.

From this discussion and the definition of the classical conditional 
expectation, we extract the following definition directly in terms of the 
quantum probability space.

\begin{definition}[Conditional expectation, finite-dimensional case]
Let $(\mathscr{N},\mathbb{P})$ be a finite-dimensional quantum probability 
space and let $\mathscr{A}\subset\mathscr{N}$ be a commutative 
$*$-subalgebra.  Then $\mathbb{P}(\cdot|\mathscr{A}):\mathscr{A}'\to
\mathscr{A}$ is called (a version of) the conditional expectation from 
$\mathscr{A}'$ onto $\mathscr{A}$ if $\mathbb{P}(\mathbb{P}(B|\mathscr{A})A)=
\mathbb{P}(BA)$ for all $A\in\mathscr{A}$, $B\in\mathscr{A}'$.
\end{definition}

As we will see in \S\ref{sec QP}, the discussion above generalizes
directly to the infinite-dimensional case.  In finite dimensions it is
convenient to give an explicit expression for the conditional expectation.
Note that a finite-dimensional $*$-algebra is a finite-dimensional linear
space.  Then $\langle A,B\rangle_{\mathbb{P}}=\mathbb{P}(A^*B)$ turns the
algebra into a pre-Hilbert space, i.e.\ it is a Hilbert space except that
$A\mapsto\langle A,A\rangle_\mathbb{P}=\|A\|^2_\mathbb{P}$ may have a
nontrivial null space. In particular, the fundamental property
$\mathbb{P}(\mathbb{P}(B|\mathscr{A})A)=\mathbb{P}(BA)$ for all
$A\in\mathscr{A}$ is precisely that of orthogonal projection from
$\mathscr{A}'$ onto the linear subspace $\mathscr{A}$, which in a
pre-Hilbert space is uniquely determined up to an event of zero
probability.  Note that the classical characterization of
$\mathbb{P}(B|\mathscr{A})$ as the least-mean-square estimate of $B$ in
$\mathscr{A}$ follows immediately.  We will elaborate on this point in
\S\ref{sec QP}. 

An explicit expression for $\mathbb{P}(B|\mathscr{A})$ is easily obtained 
if we find an orthogonal basis for $\mathscr{A}$.  Any commutative 
$*$-algebra in finite dimensions is spanned by a set of projections that 
resolve the identity. This is easily seen: in $n$ dimensions 
any self-adjoint operator is a linear combination of at most $n$ projections
that resolve the identity, and as all the operators in the $*$-algebra 
commute they must be expressible as linear combinations of the same 
projections.  Let $\mathscr{A}={\rm span}\{P_a\}$ for some set of 
projections $P_a$.  Then a version of the conditional expectation is given by
\be
	\mathbb{P}(B|\mathscr{A})=
	\sum_{P\in\{P_a\}:\mathbb{P}(P)\ne 0}
	\frac{P}{\|P\|_\mathbb{P}}
	\left\langle
	\frac{P}{\|P\|_\mathbb{P}},
	B
	\right\rangle_\mathbb{P}
	=
	\sum_{P\in\{P_a\}:\mathbb{P}(P)\ne 0}
	\frac{\mathbb{P}(PB)}{\mathbb{P}(P)}\,P.
	\label{eg cexp fd 2}
\ee 
Note what could happen if we naively fill in some $B\not\in\mathscr{A}'$. 
Then $\langle P,B\rangle_\mathbb{P}\ne\langle B,P\rangle_\mathbb{P}$ for 
some $P\in\{P_a\}$, which implies that we obtain complex coefficients in 
the sum even if $B$ is an observable.  Hence the expression does not make 
sense unless $B\in\mathscr{A}'$. 

\begin{example}    \label{eg cexp fd}
The following example serves to illustrate conditional expectations; it 
is not meant to represent a particular physical scenario.  Consider 
$\qH={\bf C}^3$, $\mathscr{N}=M_3$ and 
$\mathbb{P}(X)=\langle\psi,X\psi\rangle$ with $\psi=(1~1~1)^T/\sqrt{3}$.
Define $A,B\in\mathscr{N}$ by
$$
	A = \left( \begin{array}{ccc} 4 & 0 & 0 
	\\  0 & 4 & 0 
	\\ 0 & 0 & 5
	\end{array} \right) =  4  \left( \begin{array}{ccc} 1 & 0 & 0 
	\\  0 & 1 & 0 
	\\ 0 & 0 & 0
	\end{array} \right)
	+ 5  \left( \begin{array}{ccc} 0 & 0 & 0 
	\\  0 & 0 & 0 
	\\ 0 & 0 & 1
	\end{array} \right),\qquad
	B = \left( \begin{array}{ccc} 0 & 1 & 0 
	\\  1 & 0 & 0 
	\\ 0 & 0 & 2
	\end{array} \right).
$$
Let $\mathscr{A}$ be the $*$-algebra generated by $A$.  Then
$$
	\mathscr{A}' = \left\{ \left(   \ba{ccc}  a & b & 0 \\ 
	c & d & 0 \\ 
	0 & 0 & x 
      \ea \right) \ : \ a, b, c, d, x \in \C \right\}.
$$
Note that $\mathscr{A}'$ is not a commutative algebra, despite that every 
element of $\mathscr{A}'$ commutes with every element of $\mathscr{A}$.
As $B\in\mathscr{A}'$, we can use \er{eg cexp fd 2} to calculate
$$
	\mathbb{P}(B|\mathscr{A}) = 
	\left( \begin{array}{ccc} 1 & 0 & 0 
	\\  0 & 1 & 0 
	\\ 0 & 0 & 2
	\end{array} \right) =  1  \left( \begin{array}{ccc} 1 & 0 & 0 
	\\  0 & 1 & 0 
	\\ 0 & 0 & 0
	\end{array} \right)
	+ 2  \left( \begin{array}{ccc} 0 & 0 & 0 
	\\  0 & 0 & 0 
	\\ 0 & 0 & 1
	\end{array} \right)  \in \mathscr{A}.
$$
The observable $\mathbb{P}(B|\mathscr{A})$ is the orthogonal projection
of $B$ onto $\mathscr{A}$ with respect to the inner product $\langle 
A,B\rangle_{\mathbb{P}}=\mathbb{P}(A^*B)$.  By the projection theorem, 
$\mathbb{P}(B|\mathscr{A})$ is an element of $\mathscr{A}$ that minimizes 
the mean square error $\|B-\mathbb{P}(B|\mathscr{A})\|_\mathbb{P}$.
\qed
\end{example}

We now proceed to develop a simple probe model that reproduces the
projection postulate.  Recall that the conditional {\em probability} of an
event $P$ given a commuting event $Q$ is simply given by
$\mathbb{P}(PQ)/\mathbb{P}(Q)$.  This is equivalent to ${\bf P}(A\cap
B)/{\bf P}(B)$ by the spectral theorem, where $A$ and $B$ are the sets
corresponding to $P$ and $Q$. 

\begin{example}(Simple probe model).
We will work in a generic $n$-dimensional setting, $n<\infty$. Let
$\qH={\bf C}^n$, $\mathscr{N}=M_n$ (the set of $n\times n$ complex
matrices), and let $\mathbb{P}(X)={\rm Tr}[\rho X]$ be some state on 
$\mathscr{N}$.  Let $A,B$ be two observables in $\mathscr{N}$ that do not 
commute.  Hence we cannot measure $A$ and $B$ directly in the same 
realization.  However, we can have the system interact with an external 
probe system, in such a way that the observable $A$ is copied to some 
probe observable $A'$ after the interaction.  If $A'$ commutes with $B$, 
we interpret this procedure (like in the Stern-Gerlach example) as an 
(indirect) measurement of $A$ followed by a (direct) measurement of $B$. 

The strategy is simple.  First, we describe the probe system by a 
separate probe quantum probability space $(\mathscr{N}_p,\mathbb{P}_p)$ and 
form the composite space $(\mathscr{N}\otimes\mathscr{N}_p,
\mathbb{P}\otimes\mathbb{P}_p)$.  Next, we introduce an interaction.
Recall from \S\ref{sec BG QM} that the evolution of an isolated system is 
described by a unitary transformation.  Hence, we will choose a probe
observable $I\otimes A'$ and construct a suitable unitary operator $U$
so that the probe observable $U^*(I\otimes A')U$ after the interaction 
gives the same outcome as $A\otimes I$ would have before the interaction.
Note that by construction, the system observable $B\otimes I$ commutes 
with $I\otimes A'$ after the interaction, $[U^*(I\otimes A')U,
U^*(B\otimes I)U]=0$.  Hence we can measure them within the same realization.

We now fill out the details of this model.  Let $A=\sum_{a\in{\rm spec}(A)}
a\,P_a$, and we denote by $m$ the number of elements in ${\rm spec}(A)$ 
(the number of possible measurement outcomes).  For the probe algebra, 
we choose $\qH_p={\bf C}^m$, $\mathscr{N}_p=M_m$.  Now fix an observable 
$A'\in\mathscr{N}_p$ that has $m$ distinct measurement outcomes. Note 
that $A'=\sum_{a\in{\rm spec}(A')}a\,P_{a}'$ and that $P_{a}'$ are 
projections onto one-dimensional subspaces of $\qH_p$; hence we can fix 
an orthonormal basis of vectors $\psi_a\in\qH_p$ such that 
$P_a'=\psi_a\psi_a^*$.  Now define the operator $X_{ab}'=\psi_b\psi_a^*+
\psi_a\psi_b^*+\sum_{c\ne a,b}\psi_c\psi_c^*\in\mathscr{N}_p$ for $a\ne 
b$, and $X_{aa}'=I$; these operators switch the events $P_a'$ and $P_b'$ 
in the sense $X_{ab}'P_a'X_{ab}'=P_b'$, $X_{ab}'P_b'X_{ab}'=P_a'$, and
$X_{ab}'P_c'X_{ab}'=P_c'$ for $c\ne a,b$.  Finally, set $\mathbb{P}_p(X)=
{\rm Tr}[P_p'\rho]$ where we have fixed some $p\in{\rm spec}(A')$ at the 
outset.

Now consider the operator $U\in\mathscr{N}\otimes\mathscr{N}_p$ defined by
$U=\sum_{a\in{\rm spec}(A)}P_a\otimes X'_{ap}$.
As $(X_{ap}')^2=I$ it follows that $U^*U=UU^*=U^2=I$, i.e.\ $U$ is 
unitary.  Note that $U^*(I\otimes P_c')U=P_c\otimes P_p'+(1-P_c)\otimes 
P_c'$ if $c\ne p$, $U^*(I\otimes P_p')U=\sum_aP_a\otimes P_a'$.
We calculate $(\mathbb{P}\otimes\mathbb{P}_p)(U^*(I\otimes P_c')U(P_c\otimes 
I))/(\mathbb{P}\otimes\mathbb{P}_p)(P_c\otimes I)=1$ for every $c$, i.e., 
the conditional probability that $U^*(I\otimes A')U$ gives the outcome 
$c$, given that we have observed $A\otimes I$ with outcome $c$, is one.  
Thus the unitary interaction $U$ precisely copies the system observable 
$A$ onto the probe observable $A'$.

We can now measure the system observable $B$ after interaction with 
the probe.  In particular, let us calculate the expectation of $B$ 
conditioned on the probe measurement.  Define $\mathscr{A}$ as the 
commutative $*$-algebra generated by $U^*(I\otimes A')U$, and note that 
$U^*(B\otimes I)U\in\mathscr{A}'$.  Thus we can use \er{eg cexp fd 2} to 
calculate 
\begin{multline*}
	(\mathbb{P}\otimes\mathbb{P}_p)(U^*(B\otimes I)U|\mathscr{A})=
	\sum_{c}\frac{
		(\mathbb{P}\otimes\mathbb{P}_p)(
		U^*(B\otimes P_c')U)
	}{(\mathbb{P}\otimes\mathbb{P}_p)(U^*(I\otimes P_c')U)}\,
	U^*(I\otimes P_c')U \\
	=\sum_c
	\frac{\mathbb{P}(P_cBP_c)}{\mathbb{P}(P_c)}\,
	U^*(I\otimes P_c')U=
	\sum_c{\rm Tr}[\rho_cB]\,U^*(I\otimes P_c')U,
\end{multline*}
where $\rho_c=P_c\rho P_c/{\rm Tr}[\rho P_c]$.  This is precisely the 
projection postulate of \S\ref{sec BG QM}.

This example may be somewhat bewildering, and we encourage the reader to 
work through the procedure for a particular model (e.g.\ that of Example 
\ref{eg cexp fd}), paying particular attention to which operators do and 
do not commute.  The reader should convince himself that different 
answers are obtained if one first measures $B$, then $A$.

Finally, we note that though we have here measured $A$ through a probe 
and $B$ directly, there is no reason to stop here.  If, in addition to 
$A$ and $B$ we want to measure an observable $C$ that does not commute 
with $B$, we would introduce a second probe to measure $B$ as well.
Now suppose that $C=A$.  If we first measure $A$ through the probe, then 
measure $A$ again we would (obviously) obtain the same outcome.  However, 
if we first probe $A$, then probe $B$, and then measure $A$, we obtain a 
different outcome than that of the first measurement of $A$!  The reader 
is encouraged to work out also this case.  The reason for this phenomenon 
is that the interaction with the probe that is used for the observation 
of $B$ disturbs the system in such a way that its value of $A$ is changed.
This effect is known as ``measurement back action''.
\qed
\end{example}

The previous example, in particular the construction of the probe and the
corresponding interaction, may seem rather {\em ad hoc}, and indeed we
have only chosen this rather artificial example to reproduce the
projection postulate.  This is not a shortcoming of the theory we have
outlined, however, but rather highlights the importance of including a
reasonable model of the probe in the quantum probability space.  Indeed,
most realistic measurement setups are not of this type and the projection
postulate of \S\ref{sec BG QM} cannot be used to describe such systems.
For example, we will see in \S\ref{sec QP} that the Stern-Gerlach
measurement is only approximately described by the projection postulate.
Later we will describe even more complicated optical measurements in which
we wish to condition system observables based on the observation of
stochastic processes in continuous time (the signal from a photodetector).
It is the latter, most practically useful case where we need quantum
filtering theory. 

\begin{remark}
It is important to realize that statements like the projection postulate
do not really implement the notion of conditioning; they consist of a 
pure conditioning component and of a particular physical probe model which 
has no statistical significance.  One also finds in the literature 
generalizations of the projection postulate, called instruments, which 
implement different types of probes \cite{Dav76,AH01}.  In the quantum 
probability context of this paper it is most natural to separate the two 
parts; we will take existing probe models from physics, and concentrate 
on the calculation of the associated conditional expectations (filtering).
\qed
\end{remark}


\section{Noncommutative probability theory}
\label{sec QP}

In the finite-dimensional case, we have seen in \S\ref{sec BG} that
quantum mechanics can be modeled as a noncommutative probability theory.
In this section we present a general formulation for quantum probability
that has wide applicability. We give a general definition of quantum
probability space, prove the existence and uniqueness of conditional
expectations, and prove a quantum version of Bayes' rule that is very
helpful for quantum filtering. 

Almost all of the features of the full theory can already be seen in the
finite-dimensional case discussed in \S\ref{sec BG}; the main difficulties
in the general case are the technicalities involved in the theory of
infinite-dimensional Hilbert spaces.  This parallells the difficulties in
classical probability theory---though finite-state random variables can be
treated by almost trivial (counting, combinatoric) methods, the
description of continuous random variables requires us to upgrade our
machinery using methods of real analysis.  Similarly, the elementary
linear algebra that underlies finite-dimensional quantum probability must
be upgraded to functional analysis if we wish to treat the
infinite-dimensional case.  Conceptually, however, the two cases are very
similar, and the reader is encouraged to develop an intuitive
understanding of the finite-dimensional case before tackling the full
formalism.  For a thorough introduction to functional analysis we refer 
to the excellent textbook \cite{ReS80}.

\subsection{Quantum probability spaces}
\label{sec QP qp}

Let $\qH$ be a complex Hilbert space, and denote by $\BS{B}(\qH)$ the set
of all bounded (linear) operators on $\qH$.  We restrict ourselves (for
the time being) to bounded operators as we wish to construct $*$-algebras
of such operators: attempting to do this with unbounded operators would
get us into no end of trouble, as we would surely run into domain
problems. Recall that for $A\in\mathscr{B}(\qH)$, the usual Hilbert space
adjoint $A^*\in\BS{B}(\qH)$ is defined by
$\langle\psi,A\phi\rangle=\langle A^* \psi,\phi\rangle$
$\forall\psi,\phi\in\qH$.  With this involution $\mathscr{B}(\qH)$ is a
$*$-algebra in the sense of \S\ref{sec BG}.

We wish to introduce a structure that plays the same role as a $*$-algebra
in the finite-dimensional case.  It turns out, however, that the
$*$-algebra structure in itself is not sufficient in the
infinite-dimensional case; we need to impose an additional technical
condition in order to be able to prove an infinite-dimensional version of
the spectral Theorem \ref{thm spectral finite}.  The additional condition
has a natural interpretation which we will discuss below;  however, the
reader should not be too worried about this technicality, particularly if 
he is not familiar with nets or locally convex topologies.  In practice we
will rarely need to verify this property directly. 

\begin{definition}
	A positive linear functional $\mu:\mathscr{B}(\qH)\to{\bf C}$
	is said to be {\em normal} if $\mu(\sup_\alpha A_\alpha)=
	\sup_\alpha\mu(A_\alpha)$ for any upper bounded increasing net
	$\{A_\alpha\}$ of positive elements in $\mathscr{B}(\qH)$.
	The locally convex topology on $\mathscr{B}(\qH)$ defined by
	the family of seminorms $\{A\mapsto|\mu(A)|:\mu~{\rm normal}\}$
	is called the {\em normal topology}.
\end{definition}

For a detailed discussion of nets, locally convex topologies, etc., see 
\cite{ReS80}.

\begin{definition}[Von Neumann algebra]   \label{def vna}
	A von Neumann algebra $\mathscr{N}$ is a $*$-subalgebra
	of $\BS{B}(\qH)$ that is closed in the normal topology.
	A state $\mathbb{P}$ on $\mathscr{N}$ is normal if it is
	the restriction to $\mathscr{N}$ of a normal state on 
	$\mathscr{B}(\qH)$. 
\end{definition}

We can now extend the spectral theorem to the infinite-dimensional case,
essentially showing that commutative von Neumann algebras with normal
states are equivalent to classical probability spaces.  See e.g.\
\cite[Proposition 1.18.1]{Sak98} for a proof. 
Conceptually, we are guided by the finite-dimensional case; Theorem 
\ref{thm spectral} extends the idea of simultaneous diagonalization to 
infinite-dimensional operators. Though technically much more involved, 
the flavor of the procedure remains the same\footnote{
	The additional measure $\mu$ that shows up in the theorem has
	no direct physical significance; its job is to identify ``enough'' 
	null sets in $L^\infty(\Omega)$ so we can construct the 
	$*$-isomorphism $\iota$.  We can generally not use ${\bf P}$ for 
	this purpose as there may be projections $P\in\mathscr{C}$ with 
	$\mathbb{P}(P)=0$; if $\iota$ were to map to 
	$L^\infty(\Omega,\mathcal{F},{\bf P})$ then necessarily 
	$\iota(P)=0$ and hence $\iota$ would not be invertible.
	The precise details of the construction are never an issue,
	as we will never use $\mu$ and only prove results ${\bf P}$-a.s.
}.

\begin{theorem}[Spectral theorem]\label{thm spectral}
	Let ${\mathscr C}$ be a commutative von Neumann algebra. Then 
	there is a measure space $(\Omega, \cF, \mu)$ and a 
	$*$-isomorphism $\iota$ from $\mathscr{C}$ to $L^\infty(\Omega,
	\cF,\mu)$, the algebra of bounded measurable complex functions 
	on $\Omega$ up to $\mu$-a.s.\ equivalence.  Moreover, a
	normal state $\mathbb{P}$ on ${\mathscr C}$ defines a 
	probability measure ${\bf P}$, which is absolutely continuous with 
	respect to $\mu$, such that $\mathbb{P}(C)=E_{\bf P}(\iota(C))$ 
	for all $C\in\mathscr{C}$.
\end{theorem}

Before we continue, let us demonstrate the significance of the additional 
technical conditions on a von Neumann algebra.  First, we give an example 
of a $*$-subalgebra of $\mathscr{B}(\qH)$ that is not a von Neumann algebra.

\begin{example}   \label{eg no projections}
Let $\qH=L^2([0,1])$ and $\mathscr{A}=C([0,1])$, the commutative algebra
of continuous functions on the unit interval.  We can consider
$A\in\mathscr{A}$ as an operator on $\qH$ under pointwise multiplication,
i.e.\ $(A\psi)(x)=A(x)\psi(x)$ for every $\psi\in\qH$. Then $\mathscr{A}$
satisfies all the requirements of a von Neumann algebra except that it is
not closed in the normal topology. Indeed, one can construct, for
example, an increasing sequence of continuous functions that converges to
$\chi_{[0,1/2]}$, which is discontinuous. 

The problem is that the only indicator functions in $\mathscr{A}$ are
$\chi_\varnothing$ and $\chi_{[0,1]}$: all other indicator functions on
$[0,1]$ are discontinuous. Hence from a probabilistic point of view
$\mathscr{A}$ defines a trivial theory, as the only events in
$\mathscr{A}$ are the trivial ones.  Nonetheless $\mathscr{A}$ is much
larger than the algebra $\C$ that is generated by $\chi_\varnothing$ and
$\chi_{[0,1]}$.  Hence $\mathscr{A}$ cannot be $*$-isomorphic to the set
of measurable functions on some measure space.  The role of normal closure 
is to avoid this complication.  Indeed, this property guarantees that any 
von Neumann algebra is generated by its projections \cite{KaR83}.  \qed
\end{example}

Like normal closure, normality of the state is also required in order for
the spectral theorem to hold.  Note that for normal states the expectation
of an increasing set of observables converges to the expectation of their
least upper bound, i.e., the monotone convergence property holds.  This
corresponds to the more basic property of countable additivity.  In the
following example we construct a state which is not normal.

\begin{example}   \label{eg not normal}
Let $\qH=\ell^2(\N)$ and $\mathscr{A}=\ell^\infty(\N)$, acting on $\qH$ by 
pointwise multiplication.  $\mathscr{A}$ is closed in the normal topology, 
i.e.\ it is a commutative von Neumann algebra.  Now introduce a 
state on $\mathscr{A}$ which is given by the expression\footnote{
	Eq.\ (\ref{eq:FAstate}) does not by itself define a state,
	as there are many $A\in\mathscr{A}$ for which the limit does not 
	exist.  However, note that Eq.\ (\ref{eq:FAstate}) is well defined 
	on a linear subspace, e.g.\ $\mathscr{D}=\{A\in\mathscr{A}:
	\exists c\in\C \mbox{ s.t.\ } \lim_{n\to\infty}A(n)=c\}$.  Now 
	$\mathbb{P}$ can be extended from $\mathscr{D}$ to $\mathscr{A}$ 
	using the Hahn-Banach theorem.
}
\begin{equation}\label{eq:FAstate}
	\mathbb{P}(A)
	=\lim_{N\to\infty}\frac{1}{N}\sum_{n=1}^N A(n),
	\qquad A\in\mathscr{D}\subset\mathscr{A}
\end{equation} 
on a suitably chosen linear subspace $\mathscr{D}$.
$\mathbb{P}$ is not a normal state; to see this, let us introduce the
events $P_n\in\mathscr{A}$ defined by $(P_n\psi)(k)=\psi(k)$ if $k\le n$,
and zero otherwise.  $\{P_n\}$ is an increasing sequence of projections in 
$\mathscr{A}$ whose least upper bound is the identity $P_\infty=I$.
However, straightforward calculation shows that $\mathbb{P}(P_n)=0$ for
any finite $n$, whereas $\mathbb{P}(I)=1$.  We conclude that the state
$\mathbb{P}$ is not normal.

Note that what we have constructed is precisely the classical model of a
uniform distribution over the natural numbers $\N$.  This does not
give rise to a well-defined probability model in the sense of Kolmogorov,
however, as the uniform distribution on $\N$ does not obey the
property that the probability of a countable union of disjoint events is 
the sum of the probabilities of these events (which is exactly what went 
wrong above).  Requiring that the state be normal is equivalent to 
requiring that it gives rise to a countably additive measure \cite{KaR86}, 
which rules out our example. 
\qed
\end{example}

\begin{remark}
	Def.\ \ref{def vna} is one of many equivalent definitions
	of a von Neumann algebra.  We have emphasized normality as it
	is close to the probabilistic notion of monotone convergence.
	Normal closure turns out to be equivalent to closure in several 
	other topologies, notably the weak and strong operator topologies 
	on $\mathscr{B}(\qH)$. We will not concern ourselves with 
	topological issues in this article; see e.g.\ \cite[sec.\ 2.4]{BrR87}.
\end{remark}

The following definition should come as no surprise.

\begin{definition}[Quantum probability space]   \label{def qp space}
	A quantum probability space is a pair $({\mathscr N},\mathbb{P})$, 
	where ${\mathscr N}$ is a von Neumann algebra and $\mathbb{P}$ is 
	a normal state.
\end{definition}

The structure has precisely the same interpretation as in \S{2}, of which 
we briefly remind the reader.  In each realization we must choose a 
commutative von Neumann subalgebra $\mathscr{A}\subset\mathscr{N}$ which 
fixes the observations.  Every statistic that pertains to these 
observations is then described by the classical probability model obtained 
by applying the spectral theorem to $(\mathscr{A},\mathbb{P})$.
The equivalence between commutative quantum probability spaces and 
classical probability spaces is the foundation of the theory; a 
commutative quantum probability model {\it is} a classical probabilistic 
model, and we will often implicitly identify these two pictures.

In this article we will only use three types of von Neumann algebras.
We list these below; they will be used throughout without comment.

(i) $\mathscr{A}=\BS{B}(\qH)$ is a von Neumann algebra.  Moreover, any 
vector state on $\mathscr{A}$ (${\mathbb P}(A)=\la \psi, A\psi \ra$ for 
fixed $\psi \in \qH$), or any convex combination of vector states, is a 
normal state.  Many models from quantum mechanics are described by 
such a model.

(ii) $\mathscr{A}=L^\infty(\Omega,\cF,{\bf P})$, acting on
$\qH=L^2(\Omega,\cF,{\bf P})$ by pointwise multiplication, is a
commutative von Neumann algebra.  Moreover, any state of the form
$\mathbb{P}(X)=E_{\bf P}(X)$ is a normal state.  This is a classical
probability model.

(iii) Given $\mathscr{S}\subset\BS{B}(\qH)$, recall that $\mathscr{S}'=
\{X \in \BS{B}(\qH):XS = SX,\ \forall S\in\mathscr{S}\}$ is called the 
\emph{commutant} of $\mathscr{S}$ in $\BS{B}(\qH)$. The following theorem 
(see \cite[Theorem 5.3.1]{KaR83} for a proof) allows us to construct von 
Neumann subalgebras of $\mathscr{B}(\qH)$.

\begin{theorem}[Double commutant theorem]
	Let $\mathscr{S}\subset\BS{B}(\qH)$ be any self-adjoint set, 
	i.e.\ $S\in\mathscr{S}\Rightarrow S^*\in\mathscr{S}$.  Then 
	$\mathscr{A}=\mathscr{S}''$ is the smallest von Neumann 
	subalgebra of $\BS{B}(\qH)$ that contains $\mathscr{S}$.  In 
	particular, $\mathscr{S}$ is a von Neumann algebra iff 
	$\mathscr{S}=\mathscr{S}''$.
\end{theorem}

Given any $\mathscr{S}\subset\BS{B}(\qH)$, we call ${\rm vN}(\mathscr{S})=
(\mathscr{S}\cup\mathscr{S}^*)''$ the von Neumann algebra generated
by $\mathscr{S}$.  We will repeatedly use this construction in the
following.  For example, suppose that we decide to measure in one
realization some commuting set of observables $A_1,\ldots,A_n$.  Then
$\mathscr{A}={\rm vN}(A_1,\ldots,A_n)$ is a commutative von Neumann 
algebra which, through the spectral theorem, describes the associated 
classical probability model.  $\mathscr{A}$ is the quantum probability 
equivalent of the $\sigma$-algebra generated by a set of random variables.

\subsection{Random variables}
\label{sec QP rv}

Now that we have a general definition of a quantum probability space, 
we can develop some tools to deal with random variables.  Recall from 
\S\ref{sec BG} that any self-adjoint element of a quantum probability 
space can be decomposed into events using Eq.\ (\ref{spectral-1}), which 
gives its interpretation as an observable (random variable).  Let us show 
how to do this in the infinite-dimensional case.

Let $(\mathscr{N},\mathbb{P})$ be a quantum probability space and consider
an element $A\in\mathscr{N}$ which is self-adjoint $A=A^*$.  Then
$\mathscr{A}={\rm vN}(A)\subset\mathscr{N}$ is a commutative von Neumann
algebra.  By the spectral theorem, there is a probability space
$(\Omega,\mathcal{F},{\bf P})$ and a $*$-isomorphism $\iota$ that maps $A$
to some (measurable) random variable $a:\Omega\to{\bf R}$.  We can now do
classical probability theory; in particular, for any Borel set
$B\in\mathcal{B}$ we have the event $[a\in B]=\{\omega\in\Omega:
a(\omega)\in B\}=a^{-1}(B)\in\mathcal{F}$. To map this event back to
$\mathscr{A}$ we simply invert $\iota$; the projection corresponding to
$[a\in B]$ is denoted by $P_A(B)=\iota^{-1}(\chi_{[a\in B]})$, and we call
the map $P_A$ from $\mathcal{B}$ to the projections in $\mathscr{N}$ the
{\it spectral measure} of $A$.  But this object is a familiar one from 
functional analysis \cite{ReS80}; in fact, it is well known that we can 
express $A$ in terms of its spectral measure by
\begin{equation} \label{eq:specmeasintg}
	A = \int_{\bf R} \lambda\,P_A(d\lambda)
\end{equation}
where the integral is defined in a suitable sense \cite{ReS80}.  
Eq.\ (\ref{eq:specmeasintg}) is precisely the infinite-dimensional 
counterpart of Eq.\ (\ref{spectral-1}).  We emphasize the physical 
interpretation of $P_A(B)$: it is the event [$A$ takes a value in $B$],
which occurs with probability $\mathbb{P}(P_A(B))$.

This would be all there is to it, were it not for the fact that our
algebras contain only bounded operators (recall that unbounded operators
cannot be defined on the entire Hilbert space, and hence cannot be added
or multiplied at will). Evidently we didn't lose much by this choice, as
the probabilistic model is already contained in an algebra of bounded
operators by the spectral theorem.  An unfortunate side effect, however,
is that self-adjoint operators in the algebra can only represent {\it
bounded} random variables, whereas many observations of interest are quite
naturally unbounded (think of a Gaussian random variable).  This means
that we need to deal with unbounded observables separately.  We briefly
discuss one way of doing this.

Consider a von Neumann algebra $\mathscr{N}\subset\mathscr{B}(\qH)$.  In
general, an observable is defined by a (not necessarily bounded)
self-adjoint operator $A$ on some dense domain in $\qH$.  We need to
relate the unbounded operator $A$ to $\mathscr{N}$.  The trick we use is
remarkably simple: we compute a bounded function of $A$.  Define
$T_A=(A+iI)^{-1}$.  By elementary spectral theory \cite{ReS80}, any
self-adjoint $A$ has a real spectrum, and hence $A+iI$ is invertible with
bounded inverse.  We say that $A$ is {\it affiliated} to $\mathscr{N}$ if
$T_A\in\mathscr{N}$.  This is the equivalent of the classical notion of a
random variable that is measurable with respect to some $\sigma$-algebra
$\mathcal{G}$.  Note that every self-adjoint $A$ is affiliated to
$\mathscr{B}(\qH)$, and if $A$ is also bounded then $A$ is affiliated to
$\mathscr{N}$ iff $A\in \mathscr{N}$.

We wish to represent $A$ as a classical (unbounded) random variable.  To
this end, define the von Neumann algebra generated by $A$ as ${\rm
vN}(A)={\rm vN}(T_A)$.  Now note that $T_A$ commutes with its adjoint,
hence ${\rm vN}(A)$ is a commutative von Neumann algebra to which we can
apply the spectral theorem.  All we need to do is to ``package'' $A$ into
$T_A$, apply $\iota$, and ``unpack'' it on the other end; in other words,
we define $\iota(A) = \iota(T_A)^{-1} - i$.  Once we have done this, we
can define a spectral measure $P_A$ for $A$ in the usual way, and indeed
Eq.\ (\ref{eq:specmeasintg}) still holds even for unbounded $A$
\cite{ReS80}.  We remark that $A$ being affiliated to $\mathscr{N}$
corresponds to the fact that $P_A(B)\in\mathscr{N}$ for every
$B\in\mathcal{B}$; this is precisely the classical notion of
measurability.

Unbounded operators are a nuisance, but unfortunately they are a fact of
life in mathematical physics.  In this article, particularly in the later
sections, we will occasionally add and multiply unbounded operators
without justification; a detailed analysis of the operator domains is
beyond our scope.  Though this does not often cause trouble, the reader
should keep in mind that a fully rigorous treatment must verify that any
addition or multiplication of unbounded operators is indeed well defined.  
We quote one useful result: operators affiliated to a {\it commutative}
von Neumann algebra can be added and multiplied at will \cite[Theorem
5.6.15]{KaR83}, \cite{Nelson}.

\begin{example}  \label{eg q and p}
We take $\qH=L^2(\R)$ and $\mathscr{N}=\mathscr{B}(\qH)$.  The vector
$$
	\psi\in \qH,\qquad
	\psi(x)=(2\pi)^{-1/4}\sigma^{-1/2}\exp\left(
		-\frac{(x-\mu)^2}{4\sigma^2}
	\right)
$$
defines the (pure) state $\mathbb{P}(X)=\langle\psi,X\psi\rangle$.  Now 
consider the self-adjoint operators
$$
	(Q\psi)(x)=x\psi(x),\qquad\qquad 
	(P\psi)(x)=-i\hbar\frac{d}{dx}\psi(x),
$$
which are prototypical observables for the position $Q$ and momentum $P$ 
of a quantum particle.  Both are unbounded observables, but their domains 
include at least the set of smooth functions with compact support which 
is dense in $L^2(\R)$.

What random variables do these represent?  We can read off from the 
definition that $Q$ is a Gaussian random variable with mean $\mu$ and 
variance $\sigma^2$---as $Q$ is already in ``diagonal'' form ($Q$
is affiliated to $L^\infty(\R)\subset\mathscr{N}$), its spectral 
measure is given by
$$
	(P_Q(B)\psi)(x)=\chi_B(x)\psi(x)
$$
and it is evident that $\mathbb{P}(P_Q(B))$ is a Gaussian measure with 
mean $\mu$ and variance $\sigma^2$.  Alternatively, consider the 
characteristic function $q(k)=\mathbb{P}(e^{ikQ})$ of $Q$.  Unlike $Q$, 
$e^{ikQ}$ is a bounded operator and we can directly compute
$$
	q(k)=\langle\psi,e^{ikQ}\psi\rangle=
	(2\pi)^{-1/2}\sigma^{-1}
	\int_{-\infty}^\infty
		e^{ikx}\,
		e^{-(x-\mu)^2/2\sigma^2}\, dx=
	e^{ik\mu-k^2\sigma^2/2}
$$
which is the characteristic function of a Gaussian random variable with 
mean $\mu$ and variance $\sigma^2$.  Similarly, $e^{ikP}$ is a bounded 
operator, and we compute
$$
	p(k)=\mathbb{P}(e^{ikP})=\langle\psi,e^{ikP}\psi\rangle=
	\int_{-\infty}^\infty \psi(x)\,\psi(x+\hbar k)\,dx=
	e^{-\hbar^2k^2/8\sigma^2}
$$
which is the characteristic function of a Gaussian random variable with 
mean zero and variance $\hbar^2/4\sigma^2$.  Thus both $Q$ and $P$ are 
Gaussian random variables, but their joint distribution is undefined as 
they do not commute.  Note that we cannot choose $\sigma$ so that both $Q$ 
and $P$ have arbitrarily small variance: this is a manifestation of the 
Heisenberg uncertainty relation (compare Eq.\ \er{hu-1}).   \qed
\end{example}

The following example plays a central role in the physics of harmonic
oscillators; we will encounter a very similar construction later for
continuous-time quantum stochastic processes.  We will need the following
classic result (see e.g.\ \cite{KaR83} for a proof).

\begin{theorem}[Stone's theorem]\label{thm stone}
	Let $\mathscr{N}$ be a von Neumann algebra and let
	$\{U_t\}_{t\in{\bf R}}\subset\mathscr{N}$ be a group of unitary 
	operators that is strongly continuous. Then there is a unique
	self-adjoint $A$ affiliated to $\mathscr{N}$, the {\em Stone 
	generator}, such that $U_t=e^{itA}$.
\end{theorem}

\begin{example} \label{ex pos}
Let $\qH=\ell^2({\bf N})$ and $\mathscr{N}=\mathscr{B}(\qH)$. Define the
complete orthonormal basis $\{\psi_n,~n=0,1,\ldots\}\subset\qH$, where
$\psi_n(k)=1$ if $k=n$ and $\psi_n(k)=0$ otherwise.  Moreover, we define
for every $\alpha\in{\bf C}$ the {\em exponential vector}
$e(\alpha)\in\qH$ by $e(\alpha)(k)=\alpha^k/\sqrt{k!}$, and we remark that
the linear span $\mathsf{D}$ of all exponential vectors is dense in $\qH$.  
The normalized exponential vectors $e(\alpha)e^{-|\alpha|^2/2}$ are called
coherent vectors, and can be used to define the coherent states
$\mathbb{P}_\alpha(X)=\langle e(\alpha),Xe(\alpha)\rangle\,e^{-|\alpha|^2}$.

The simplest random variable we can investigate is defined by
$(\lambda\psi)(k)=k\psi(k)$---i.e.\ this is the natural diagonal operator
affiliated to $\ell^\infty({\bf N})\subset\mathcal{N}$.  The spectral
measure of $\lambda$ is given by $(P_\lambda(B)\psi)(k)=\chi_B(k)\psi(k)$,
from which we obtain directly
$$
        \mathbb{P}_\alpha(P_\lambda(B))=
        \langle e(\alpha),P_\lambda(B)e(\alpha)\rangle
                \,e^{-|\alpha|^2}=
        \sum_{k\in B}\frac{e^{-|\alpha|^2}(|\alpha|^2)^k}{k!}.
$$
Thus evidently, $\lambda$ is a Poisson-distributed random variable with
intensity $|\alpha|^2$.

Can we find other interesting observables affiliated to $\mathscr{N}$?
In many cases, physically relevant observables are found to be the Stone
generators of particular unitary symmetry groups; see e.g.\ \cite{Hol82}
for a lucid discussion.  Let us try to implement this procedure with the
two-dimensional translation group.  As a first attempt, let us define a
translation operator by $D_\gamma e(\alpha)=e(\alpha+\gamma)\,
e^{|\alpha|^2/2-|\alpha+\gamma|^2/2}$ for $\gamma\in{\bf C}$; the
constant factor ensures that $\|D_\gamma e(\alpha)\|=\|e(\alpha)\|$, as
must be the case for any unitary operator.  Unfortunately, $D_\gamma$ is
not in fact unitary; a straightforward calculation shows
$$
        \langle e(\beta),D_\gamma^*D_\gamma e(\alpha)\rangle=
        \langle D_\gamma e(\beta),D_\gamma e(\alpha)\rangle=
        e^{\beta^*\alpha}
        e^{i\,{\rm Im}(\beta^*\gamma)-i\,{\rm Im}(\alpha^*\gamma)}
$$
which contradicts unitarity $D_\gamma^*D_\gamma=I$, i.e.\ $\langle
e(\beta),D_\gamma^*D_\gamma e(\alpha)\rangle=\langle e(\beta),
e(\alpha)\rangle=e^{\beta^*\alpha}$.  To fix this, define the
{\em Weyl operator}
$$
        W_\gamma e(\alpha)=e(\alpha+\gamma)\,
        e^{|\alpha|^2/2-|\alpha+\gamma|^2/2}e^{i\,{\rm Im}(\alpha^*\gamma)}=
        e(\alpha+\gamma)\,e^{-\gamma^*\alpha-|\gamma|^2/2}.
$$
The Weyl operator is unitary, and provides a projective unitary
representation \cite{Hol82} in the sense that $W_\alpha
W_\beta=W_{\alpha+\beta}e^{i\,{\rm Im}(\beta^*\alpha)}$.  Note that it is
sufficient to define the action of $W_\alpha$ only on exponential vectors;
we can then extend to $\mathsf{D}$ by linearity, and as $\mathsf{D}$ is
dense and $W_\alpha$ is bounded the Weyl operators are uniquely extended
to all of $\qH$.

Now fix $\beta\in{\bf C}$ and consider the unitary group 
$\{W_{t\beta}\}_{t\in{\bf R}}$.  This group is continuous
($W_{t\beta}e(\gamma)\to e(\gamma)$ as $t\to 0$) and hence by Stone's
theorem, there exists a self-adjoint operator $B_\beta$ such that
$W_{t\beta}=e^{itB_\beta}$.  Finding the distribution of the observable
$B_\beta$ is straightforward, as the chararcteristic function of $B_\beta$
is given by
$$
        b_\beta(k)=\mathbb{P}_\alpha(W_{k\beta})=
        \langle e(\alpha),e(\alpha+k\beta)\rangle
                \,e^{-k\beta^*\alpha-k^2|\beta|^2/2-|\alpha|^2}=
        e^{2ik\,{\rm Im}(\alpha^*\beta)-k^2|\beta|^2/2}.
$$
Hence $B_\beta$ is a Gaussian random variable with mean $2\,{\rm
Im}(\alpha^*\beta)$ and variance $|\beta|^2$.

Our next task is to obtain an explicit representation of $B_\beta$.
We proceed as follows:
$$
        B_\beta e(\alpha)=
        \left.\frac{1}{i}\frac{d}{dt}W_{t\beta}e(\alpha)\right|_{t=0}=
        i\beta^*\alpha \,e(\alpha)-i
        \left.\frac{d}{dt}e(\alpha+t\beta)\right|_{t=0}.
$$
One can verify explicitly that this expression makes sense, i.e.\
$B_\beta e(\alpha)\in\qH$.  Note that we cannot extend $B_\beta$
to all of $\qH$, as $B_\beta$ is unbounded.  However, we see that
the domain of $B_\beta$ contains at least the exponential domain
$\mathsf{D}$.

Let us introduce the following notation.  Define $q=B_i$, $p=B_{-1}$, and
$a=(q+ip)/2$.  Note that $q$ and $p$ are self-adjoint by Stone's theorem,
whereas $a$ has the adjoint $a^*=(q-ip)/2$.  Moreover, we find that
$a\,e(\alpha)=\alpha\,e(\alpha)$.  But then
$$
        (a\,e(\alpha))(k)=\alpha\,\frac{\alpha^k}{\sqrt{k!}}=
        \sqrt{k+1}\,\frac{\alpha^{k+1}}{\sqrt{(k+1)!}} =
        \sqrt{k+1}\,e(\alpha)(k+1).
$$
This implies that we can extend the domain of $a$ to include also the
$\{\psi_n\}$ by defining $a\,\psi_{k+1}=\sqrt{k+1}\,\psi_k$ (where
$a\,\psi_0=0$). Furthermore, from
$$
        \langle\psi_m,a^*\psi_k\rangle=\langle a\psi_m,\psi_k\rangle=
        \sqrt{m}\,\delta_{(m-1)k}=\sqrt{k+1}\,\delta_{m(k+1)}
$$
we can read off $a^*\psi_k=\sqrt{k+1}\,\psi_{k+1}$.  $a^*$ is known as the
creation (or raising) operator and $a$ as the annihilation (or lowering)
operator.

Finally, note that $\lambda=a^*a$.  From a classical probability point
of view this is very remarkable indeed.  Not only do both Poisson and
Gaussian random variables emerge from the same state $\mathbb{P}_\alpha$,
but there is even a {\it continuous} map 
$q,p\mapsto(q-ip)(q+ip)/4=\lambda$
that transforms two Gaussian random variables into a Poisson random
variable.  One could never continuously transform a continuous classical
random variable into a discrete classical random variable; however, we
get away with it here because $p$, $q$ and $\lambda$ do not commute with
one another. Thus in each realization we can choose to measure either a
discrete or a continuous random variable, but not both.
\qed
\end{example}

\begin{remark}\label{remark stariso pos}
Though presented rather differently, the last two examples are in fact
$*$-isomorphic in the case that $\sigma^2=\tfrac{1}{2}$ in the first
example.  For example, if $\alpha\in{\bf R}$ we can map $p\mapsto
2^{1/2}\hbar^{-1}P$, $q\mapsto 2^{1/2}Q$, and $\mathbb{P}_\alpha\mapsto
\mathbb{P}_{\mu=2^{1/2}\alpha,\sigma=2^{-1/2}}$.  From the expression for
$b_\beta(k)$ we see that in a coherent state both $p$ and $q$ must have
the same variance.  In the first example we allowed for the variance of
$Q$ to shrink, though this necessarily increases the variance of $P$.  
This results in a ``squeezed state'' which can also be introduced in the 
context of the second example.  We will not construct such states here; in 
the following, we will only use coherent states. \qed
\end{remark}

\subsection{Conditional expectation}
\label{sec QP ce}
 
We now consider conditional expectations, following the treatment of
\cite{BH05}.  The following definition is identical to the one in
\S\ref{sec BG}.

\begin{definition}[Conditional expectation]\label{de conditional expectation}
Let $(\mathscr{N},\mathbb{P})$ be a quantum probability space and let 
$\mathscr{A}\subset\mathscr{N}$ be a commutative von Neumann subalgebra.  
Then the map $\mathbb{P}(\cdot|\mathscr{A}):\mathscr{A}'\to\mathscr{A}$ is 
called (a version of) the conditional expectation from $\mathscr{A}'$ onto 
$\mathscr{A}$ if $\mathbb{P}(\mathbb{P}(B|\mathscr{A})A)=
\mathbb{P}(BA)$ for all $A\in\mathscr{A}$, $B\in\mathscr{A}'$.
\end{definition}

We briefly recall the significance of $\mathscr{A}'$.  $\mathscr{A}$ 
is the algebra generated by our observations: it must be commutative, as 
we cannot observe incompatible events in a single experiment.  We now wish 
to find the conditional statistics of an observable $B$ that is not
affiliated to $\mathscr{A}$.  However, as we have already observed
$\mathscr{A}$, this is only sensible if $B$ commutes with every element in
$\mathscr{A}$---there would be no physical way to test our predictions if
we could not subsequently measure $B$ in the same realization.

\begin{remark}  \label{rmk ce affiliated}
Recall that if $B=B^*$ we can use the spectral theorem to obtain 
explicitly $\mathbb{P}(B|\mathscr{A})=\iota^{-1}(E_{\bf P}(\iota(B)|
\sigma\{\iota(\mathscr{A})\}))$.  This representation extends even to the 
case that $B$ is an unbounded self-adjoint operator that is affiliated to 
$\mathscr{A}'$.  For simplicity we will discuss below the properties of 
$\mathbb{P}(B|\mathscr{A})$ assuming that $B$ is bounded, but with 
suitable care the treatment extends also to the unbounded case.
\qed
\end{remark}

\begin{remark}  \label{rmk ce general}
A more general definition (see e.g.\ \cite{Tak71}), of which Definition
\ref{de conditional expectation} is a special case, is often used in
quantum probability.  Unlike our definition, which is motivated by
statistical inference and filtering, the more general ``conditional
expectation'' allows for conditioning on noncommutative algebras and does
not have a direct statistical interpretation.  The more general definition
is used e.g.\ in the theory of noncommutative Markov processes
\cite{Kum85}.  We will not dwell on this further. \qed 
\end{remark}

\begin{theorem}   \label{thm ce}
	The conditional expectation of Definition {\rm\ref{de conditional 
	expectation}} exists and is unique with probability one 
	{\rm(}any two versions $P$ and $Q$ of 
	$\mathbb{P}(B|\mathscr{A})$ satisfy $\|P-Q\|_\mathbb{P}=0$, where
	$\|X\|^2_\mathbb{P}=\mathbb{P}(X^*X)$.{\rm)}
	Moreover, $\mathbb{P}(B|\mathscr{A})$ is the least mean 
	square estimate of $B$ given $\mathscr{A}$ in the sense that
	$\|B-\mathbb{P}(B|\mathscr{A})\|_\mathbb{P}\le \|B-A\|_\mathbb{P}$ 
	for all $A\in\mathscr{A}$.
\end{theorem}

\begin{proof}~

(i) {\it Existence.} We have already established that for self-adjoint
$B\in\mathscr{A}'$, we can explicitly define a $\mathbb{P}(B|\mathscr{A})$
that satisfies the conditions of Definition \ref{de conditional expectation}
using the spectral theorem.  The classical conditional expectation exists, 
and moreover the conditional expectation of a bounded random variable is 
bounded.  Hence $\mathbb{P}(B|\mathscr{A})$ exists in $\mathscr{A}$ for 
self-adjoint $B\in\mathscr{A}'$.  But any $B\in\mathscr{A}'$ can be 
written as $B=B_1+iB_2$ with self-adjoint $B_1=(B+B^*)/2$ and 
$B_2=i(B^*-B)/2$. As $\mathbb{P}(B_1|\mathscr{A})$ and 
$\mathbb{P}(B_2|\mathscr{A})$ exist and 
$\mathbb{P}(B|\mathscr{A})=\mathbb{P}(B_1|\mathscr{A})
+i\mathbb{P}(B_2|\mathscr{A})$ satisfies the conditions of Definition
\ref{de conditional expectation}, existence is proved.

(ii) {\it Uniqueness w.p.\ one.} Define the pre-inner product
$\langle X, Y\rangle = \mathbb{P}(X^*Y)$ on ${\mathscr A}'$ (it might have
nontrivial kernel). Then $\langle A,B-\BB{P}(B|{\mathscr A})\rangle =
\mathbb{P}(A^*B)-\mathbb{P}(A^*\BB{P}(B|{\mathscr A}))=0$ for all $A \in 
{\mathscr A}$ and $B \in {\mathscr A}'$, i.e.\ $B -\BB{P}(B|{\mathscr A})$ 
is orthogonal to ${\mathscr A}$.  Now let $P$ and $Q$ be two versions of 
$\BB{P}(B|{\mathscr A})$.  It follows that $\langle A,P-Q\rangle=0$ for 
all $A\in{\mathscr A}$.  But $P-Q\in{\mathscr A}$, so $\langle P-Q,P-Q\rangle=
\|P-Q\|_\mathbb{P}^2=0$.

(iii) {\it Least squares.}  Let $P$ be a version of $\BB{P}(B|{\mathscr 
A})$.  Then for all $K\in\mathscr{A}$
\begin{equation*}
	\|B-K\|_\mathbb{P}^2 = \|B-P+P-K\|_\mathbb{P}^2 
		= \|B-P\|_\mathbb{P}^2 + \|P-K\|_\mathbb{P}^2 
	\ge \|B-P\|_\mathbb{P}^2
\end{equation*}
where, in the second step, we used that $(B-\BB{P}(B|{\mathscr A})) \perp
(\BB{P}(B|{\mathscr A})-K)\in {\mathscr A}$. 
\qquad
\end{proof}

\begin{remark}  \label{rmk ce classical}
The usual elementary properties of classical conditional expectations and 
their proofs \cite{Wil91} carry over directly.  In particular, we have 
linearity, positivity, invariance of the state 
$\mathbb{P}(\mathbb{P}(B|\mathscr{A}))=\mathbb{P}(B)$, invariance of
$\mathscr{A}$ ($\mathbb{P}(B|\mathscr{A})=B$ if $B\in\mathscr{A}$), the
tower property $\mathbb{P}(\mathbb{P}(B|\mathscr{A})|\mathscr{C})=
\mathbb{P}(B|\mathscr{C})$ if $\mathscr{C}\subset\mathscr{A}$, the module
property $\mathbb{P}(AB|\mathscr{C})=B\,\mathbb{P}(A|\mathscr{C})$ for
$B\in\mathscr{C}$, etc.  As an example, let us prove linearity.  It
suffices to show that $Z=\alpha\,\mathbb{P}(A|\mathscr{C})+
\beta\,\mathbb{P}(B|\mathscr{C})$ satisfies $\mathbb{P}(ZC)=
\mathbb{P}((\alpha A+\beta B)C)$ for all $C\in\mathscr{C}$.  But this is 
immediate from the linearity of $\mathbb{P}$ and Definition \ref{de 
conditional expectation}. 
\qed
\end{remark}

\subsection{The Bayes formula}
\label{sec QP bayes}

In \S\ref{sec BG} we were able to calculate conditional expectations
explicitly as all algebras were finite-dimensional.  In most physical 
situations, however, at least the probe (and often the system as well) 
admits continuous observables and therefore we must deal with 
infinite-dimensional algebras.  In this case it is usually not so simple
to calculate the conditional expectations directly; however, the following
Bayes-type formula will be of considerable assistance.

\begin{lemma}[Bayes formula \cite{BH05}]\label{thm KS}
Let $\mathscr{C}$ be a commutative von Neumann algebra and let 
$\mathscr{C}'$ be equipped with a normal state $\mathbb{P}$.  Choose
$V\in{\mathscr C}'$ such that $V^*V>0$ and $\BB{P}(V^*V)=1$.  Then we can 
define a new state on ${\mathscr C}'$ by $\BB{Q}(A) = \BB{P}(V^*AV)$ and
$$
	\BB{Q}(X|{\mathscr C}) = \frac{\BB{P}(V^*XV|{\mathscr C})}
	{\BB{P}(V^*V|{\mathscr C})},\qquad X\in {\mathscr C}'.
$$
\end{lemma}

\begin{proof} 
Let $K$ be an element of ${\mathscr C}$. For all $X \in {\mathscr C}'$, we 
can write
\begin{multline*}
	\BB{P}(\BB{P}(V^*XV|{\mathscr C})K)
	= \BB{P}(V^*XKV)
	= \BB{Q}(XK)
	= \BB{Q}(\BB{Q}(X|{\mathscr C})K) \\ 
   	= \BB{P}(V^*V\BB{Q}(X|{\mathscr C})K)
	= \BB{P}(\BB{P}(V^*V\BB{Q}(X|{\mathscr C})K|{\mathscr C}))
	= \BB{P}(\BB{P}(V^*V|{\mathscr C})\BB{Q}(X|{\mathscr C})K).
\end{multline*}
As this holds for all $K\in{\mathscr C}$, and as by construction the 
conditional expectations are elements of ${\mathscr C}$, we conclude that
$\|\BB{P}(V^*XV|{\mathscr C})-\BB{P}(V^*V|{\mathscr C})\BB{Q}(X|{\mathscr 
C})\|_\mathbb{P}=0$, or equivalently $\BB{P}(V^*XV|{\mathscr C})=
\BB{P}(V^*V|{\mathscr C})\BB{Q}(X|{\mathscr C})$ $\mathbb{P}$-a.s.
\qquad
\end{proof}

We now have sufficient tools to deal with the Stern-Gerlach experiment 
described in \S\ref{sec BG}.  Though the following example is not of much 
practical importance, it demonstrates the use of the Bayes theorem in a 
concrete setting.  We will use a very similar ``reference probability 
method'' to obtain filtering equations later on.

\begin{example}(Stern-Gerlach experiment).   \label{eg stern}
Consider an atom with two degrees of freedom: a spin degree of freedom 
$\mathscr{N}_\mu=\mathscr{B}({\bf C}^2)$ carrying the observables 
$\sigma_x$, $\sigma_z$ etc., and a single spatial degree of freedom 
$\mathscr{N}_x=\mathscr{B}(\ell^2({\bf N}))$ with the affiliated 
position $q$ and momentum $p$ observables defined\footnote{
	We saw in Remark \ref{remark stariso pos} that this description
	is $*$-isomorphic to the usual definition of position and momentum 
	up to some numerical constants.  These are not of essence, 
	however, as they just correspond to a change of units in which
	we measure position and momentum.  A little more care must be 
	taken if we wish to make quantitative predictions on the outcomes 
	of actual experiments; we will not worry about this, however,
	and work in arbitrary units.
}
in Example \ref{ex pos} (we use the notations of that example).
The total algebra describing the atom is then $\mathscr{N}=\mathscr{N}_\mu
\otimes\mathscr{N}_x$.  Initially the spin and position/momentum of the 
atom are uncorrelated; hence we work with the state 
$\mathbb{P}=\mathbb{P}_\mu\otimes\mathbb{P}_0$, where $\mathbb{P}_\mu$ is
an arbitrary spin state and $\mathbb{P}_0(X)=\langle\psi_0,X\psi_0\rangle=
\langle e(0),Xe(0)\rangle$.  The latter implies that initially $I\otimes 
q$ and $I\otimes p$ (which we will interpret as position and momentum in 
the $z$-direction) have zero mean and unit variance.

To measure the spin, we apply a magnetic field gradient that is linear in 
$q$ for some fixed period of time. The resulting force on the particle 
will cause its momentum to change; an observation of the momentum of the 
particle after the interaction should thus provide a measurement of its 
spin $\sigma_z$.  In other words, the atomic spatial degree of freedom 
acts as a probe for the atomic spin degree of freedom.  The action of the 
magnetic field is described by the unitary\footnote{
	This is the solution of Eq.\ (\ref{schrodinger}) at some fixed 
	time $t$ for a suitable interaction Hamiltonian $H$.
}
$$
        U=\exp(i\kappa\,\sigma_z\otimes q)=P_{z,1}\otimes e^{i\kappa q}+
                        P_{z,-1}\otimes e^{-i\kappa q}=
                P_{z,1}\otimes W_{i\kappa}+P_{z,-1}\otimes W_{-i\kappa}
$$
where $\kappa\in{\bf R}$ is the field gradient.  Let us thus begin by 
calculating the characteristic function of $U^*(I\otimes p)U$, the 
momentum of the atom after the interaction:
\begin{multline*}
        \mathbb{P}(e^{ik\,U^*(I\otimes p)U})
                =\mathbb{P}(U^*(I\otimes W_{-k})U)
                =\mathbb{P}_\mu(P_{z,1})\,
                        \mathbb{P}_x(W_{-i\kappa}W_{-k}W_{i\kappa}) \\ +
                 \mathbb{P}_\mu(P_{z,-1})\,
                        \mathbb{P}_x(W_{i\kappa}W_{-k}W_{-i\kappa})
                =\mathbb{P}_\mu(P_{z,1})\,e^{2i\kappa k-k^2/2}+
                 \mathbb{P}_\mu(P_{z,-1})\,e^{-2i\kappa k-k^2/2}.
\end{multline*}
Hence the momentum of the atom after the interaction is distributed as a
sum of two Gaussians of unit variance and means $2\kappa$ and $-2\kappa$,
which are weighted respectively by $\mathbb{P}_\mu(P_{z,1})$ and
$\mathbb{P}_\mu(P_{z,-1})$.  Note that we cannot perfectly resolve the 
spin-up and down states using a Stern-Gerlach measurement; as the tails 
of the two Gaussians overlap, there is always a nonzero probability that 
we assign the wrong spin to the atom by looking e.g.\ at the sign of the 
observed momentum.  However, the error probability becomes very small when 
the gradient $\kappa$ is large.

After the interaction, we may want to measure a spin observable
$\sigma\in\mathscr{N}_\mu$ that does not necessarily commute with
$\sigma_z$ (e.g.\ $\sigma_x$).  To describe this, let us calculate
$\mathbb{P}(U^*(\sigma\otimes I)U|{\rm vN}(U^*(I\otimes p)U))$, the
conditional expectation of the spin observable $\sigma$ after the
interaction given our observation of the momentum of the atom.

We begin by using the following elementary property: if $U$ is a unitary
operator and we define the state $\mathbb{Q}(X)=\mathbb{P}(U^*XU)$, then
$\mathbb{P}(U^*XU|U^*\mathscr{C}U)= U^*\mathbb{Q}(X|\mathscr{C})U$ (this
can be verified directly using Definition \ref{de conditional 
expectation}). Thus we obtain
$$
        \mathbb{P}(U^*(\sigma\otimes I)U|{\rm vN}(U^*(I\otimes p)U)) =
        U^*\mathbb{Q}(\sigma\otimes I|{\rm vN}(I\otimes p))U.
$$
We would like to apply the Bayes rule to $\mathbb{Q}(\sigma\otimes I|{\rm
vN}(I\otimes p))$.  As $U$ does not commute with $I\otimes p$, however,
the Bayes rule does not apply in this form.

Fortunately we can circumvent this problem using the following trick.
Using the Baker-Campbell-Hausdorff formula, we can rewrite $e^{i\kappa q}$ 
as
$$
        e^{i\kappa q}=
        e^{i\kappa(a+a^*)}=
        e^{-\kappa^2/2}e^{i\kappa a^*}e^{i\kappa a}.
$$
Beware that the Baker-Campbell-Hausdorff formula technically only holds 
for exponentials of bounded operators; thus here and below there will be 
domain issues, but these can be resolved with suitable care.  As 
$a\,\psi_0=0$, we can write
$$
        e^{i\kappa q}\psi_0
        =e^{-\kappa^2/2}e^{i\kappa a^*}e^{i\kappa a}\psi_0=
        e^{-\kappa^2/2}e^{i\kappa a^*}\psi_0=
        e^{-\kappa^2/2}e^{i\kappa a^*}e^{-i\kappa a}\psi_0=
        e^{-\kappa^2}e^{\kappa p}\psi_0.
$$
We obtain
$$
        \mathbb{P}_0(e^{-i\kappa q}Xe^{i\kappa q})
        =\langle e^{i\kappa q}\psi_0,Xe^{i\kappa q}\psi_0\rangle
        =e^{-2\kappa^2}\langle e^{\kappa p}\psi_0,Xe^{\kappa p}\psi_0\rangle
        =e^{-2\kappa^2}\mathbb{P}_0(e^{\kappa p}Xe^{\kappa p}).
$$
It follows that we can equivalently replace $U$ by $V$:
$$
        \mathbb{Q}(X)=\mathbb{P}(U^*XU)=
                \mathbb{P}(V^*XV),\quad
        V=e^{-\kappa^2}e^{\kappa\,\sigma_z\otimes p}=
        e^{-\kappa^2}(P_{z,1}\otimes e^{\kappa p}+
                P_{z,-1}\otimes e^{-\kappa p}).
$$
$V$ is not unitary, but it does commute with $I\otimes p$.
Hence the Bayes rule gives
$$
        \mathbb{P}(U^*(\sigma\otimes I)U|{\rm vN}(U^*(I\otimes p)U)) =
        \frac{U^*\mathbb{P}(V^*(\sigma\otimes I)V|{\rm vN}(I\otimes p))U}
        {U^*\mathbb{P}(V^*V|{\rm vN}(I\otimes p))U}.
$$
We can now use the module property and independence of $\sigma\otimes 
I$ and $I\otimes p$ under $\mathbb{P}$ to calculate explicitly the 
numerator and denominator; elementary manipulations give
\begin{multline*}
        \mathbb{P}[U^*(\sigma\otimes I)U|{\rm vN}(U^*(I\otimes p)U)] = \\
        \frac{
                \mathbb{P}_\mu(P_{z,1}\sigma P_{z,1})  
                        e^{2\kappa\,U^*(I\otimes p)U}
                +\mathbb{P}_\mu(P_{z,-1}\sigma P_{z,-1})  
                        e^{-2\kappa\,U^*(I\otimes p)U}
                +2\,{\rm Re}\,\mathbb{P}_\mu(P_{z,-1}\sigma P_{z,1})
        }{
                \mathbb{P}_\mu(P_{z,1})  
			e^{2\kappa\,U^*(I\otimes p)U}
                +\mathbb{P}_\mu(P_{z,-1})  
			e^{-2\kappa\,U^*(I\otimes p)U}
        }.
\end{multline*}
By definition $\mathbb{P}(U^*(\sigma\otimes I)U|{\rm vN}(U^*(I\otimes 
p)U))$ is affiliated to ${\rm vN}(U^*(I\otimes p)U)$, and indeed the
expression above is simply a function of $U^*(I\otimes p)U$.  If we
observe $U^*(I\otimes p)U$ and obtain the value $\tilde p$, then the
spectral theorem tells us that the conditional expectation takes the value
given by the expression above if we simply substitute $\tilde p$ for
$U^*(I\otimes p)U$.  Note that the formula is not equivalent to the
one given by the projection postulate for a measurement of $\sigma_z$.  
For large $\kappa$, however, we obtain approximately the projection 
postulate expression, and this becomes exact as $\kappa\to\infty$.
\qed
\end{example}


\section{Stochastic processes and quantum It\^o calculus}
\label{sec Fock}

After a general introduction to quantum probability, we now turn to one
particular quantum probability space which we will use throughout the
remainder of the article. In \S\ref{sec filtering} we shall argue that
this model appropriately describes the quantum electromagnetic field and
its interaction with matter.  In the laboratory, the electromagnetic field
can be measured by devices like photodetectors which can produce an
electric current or even a discrete photocount. The statistics of data
records from such experiments are well approximated by the model
considered here.  The model is rich and we will discover that it contains
many interesting classical stochastic processes, i.e.\ a whole family of
Poisson and Wiener processes. However, these processes do not commute with
each other. An extension of the It\^o calculus, due to Hudson and
Parthasarathy \cite{HuP84}, unites all these processes in one
noncommutative stochastic calculus.

\subsection{Poisson processes on Fock space}

The theory we are about to discuss can be approached from many sides;  
here we have chosen to get started by finding a quantum probability space
that naturally admits a Poisson process, and build the theory from there.  
As we have a particular classical process in mind, the general theory
gives a hint as to how we could proceed.  First, we define the process on
a classical space $(\Omega,\mathcal{F},{\bf P})$; equivalently, we can
form the algebra $\mathscr{A}=L^\infty(\Omega,\mathcal{F},{\bf P})$ acting
on $\qH=L^2(\Omega,\mathcal{F},{\bf P})$ by pointwise multiplication, with
a suitable state $\mathbb{P}$, and represent the process as a family of
observables affiliated to $\mathscr{A}$.  To create a noncommutative
model, we could now broaden our horizon and consider $\mathscr{N}=
(\mathscr{B}(\qH),\mathbb{P})$ rather than just $\mathscr{A}$.  Obviously
such a construction does not necessarily carry a physical interpretation;  
this must be considered separately, see \S\ref{sec filtering}.  For the
time being, however, we will use this convenient construction to provide
us with a rich quantum stochastic model.  The following discussion is
heavily inspired by the work of Maassen \cite{Maa85}.

Consider a classical Poisson process on a finite time interval $[0,T]$.
We wish to describe the space of paths $\Omega$.  This is not difficult; 
a Poisson process on a finite time interval has (a.s.)\ finitely many 
jumps $n$.  Hence we can specify every relevant path by specifying its 
jump times.  Let us thus introduce
\begin{equation}
	\Omega=\bigcup_{n=0}^\infty\Omega_n,	\qquad
	\Omega_0=\{\varnothing\},\quad
	\Omega_n=\{\{t_1,\ldots,t_n\}:t_1<t_2<\ldots<t_n\in [0,T]
		\}.
\end{equation}
In other words, $\Omega$ is the set of ordered sequences in $[0,T]$ with a
finite number of elements.  We still need to introduce a $\sigma$-algebra
$\mathcal{F}$ and a measure ${\bf P}$.  To this end, consider $\Omega_n$
as a subset of the cube $([0,T]^n,e^{-T}\mu_n)$ where $\mu_n$ is the
Lebesgue measure, so that $\Omega_n$ inherits a $\sigma$-algebra
$\mathcal{F}_n$ and a measure ${\bf P}_n$ from the cube.  Under ${\bf
P}_n$ the jump times $t_1,\ldots,t_n$ are uniformly distributed (as must
be the case for a Poisson process with fixed rate) and ${\bf
P}_n(\Omega_n)= T^ne^{-T}/n!$.  The measure ${\bf P}$ induced on $\Omega$
is precisely the probability measure of a Poisson process with unit rate.

We now introduce the Hilbert space $\mathsf{F}=L^2(\Omega,\mathcal{F},{\bf 
P})$.  It is called the {\it symmetric} or {\it Boson Fock space}, and 
plays a central role in the following.  We will also need the spaces 
$\mathsf{F}_{t]}$, $\mathsf{F}_{[t}$ and $\mathsf{F}_{[s,t]}$, defined 
identically to $\mathsf{F}$ except that the interval $[0,T]$ is replaced 
by $[0,t]$, $[t,T]$ and $[s,t]$, respectively.  It is not difficult to see 
that for any $0<s<t<T$ we have\footnote{
	A more precise statement would be something like
	$\Omega=\Omega_{s]}\times\Omega_{(s,t]}\times\Omega_{(t}$;
	however, the only paths for which this makes a difference are
	those that have jumps exactly at times $s$ or $t$, which is a set 
	of ${\bf P}$-measure zero.  For notational simplicity, we are
	free to always use closed time intervals $[s,t]$.
} $\Omega=\Omega_{s]}\times\Omega_{[s,t]}\times
\Omega_{[t}$, and as the Poisson process has independent increments the 
measure splits up similarly.  It follows that
\begin{equation}
\label{eq:guichardet}
	\mathsf{F}=\mathsf{F}_{s]}\otimes\mathsf{F}_{[s,t]}\otimes
		\mathsf{F}_{[t}\qquad\forall\,0<s<t<T.
\end{equation}
This important property is known as a continuous tensor product structure;
it will play a key role in the definition of quantum stochastic integrals, 
as it gives a natural notion of adaptedness.  Indeed, the algebra 
$\mathscr{W}=\mathscr{B}(\qF)$ splits up accordingly,
\begin{equation}
	\mathscr{W}=\mathscr{W}_{s]}\otimes\mathscr{W}_{[s,t]}\otimes
		\mathscr{W}_{[t}=\mathscr{B}(\mathsf{F}_{s]})
		\otimes\mathscr{B}(\mathsf{F}_{[s,t]})\otimes
		\mathscr{B}(\mathsf{F}_{[t}).
\end{equation}
A process of operators $\{X_t\}$ affiliated to $\mathscr{W}$ is said to be 
{\it adapted} if $X_t$ is affiliated to $\mathscr{W}_{t]}$ for every $t$; 
equivalently, $X_t$ is of the form $X_{t]}\otimes I$ as an operator on 
$\qF_{t]}\otimes\qF_{[t}$.

Next, let us introduce a set of interesting vectors.  The reader should 
keep in mind Example \ref{ex pos} which is conceptually quite similar.  
Let $f\in L^\infty([0,T])$ be a complex Lebesgue measurable function.  
Then we can define the {\it exponential vector}
\begin{equation}
\label{eq:exponential}
	e(f)(\varnothing)=1,\quad
	e(f)(\tau)=\prod_{t\in\tau}f(t),\qquad
	f\in L^\infty([0,T]).
\end{equation}
It is not difficult to verify that $e(f)\in\mathsf{F}$, as
$$
	\langle e(g),e(f)\rangle=
	\sum_{n=0}^\infty\frac{e^{-T}}{n!}\left(\int_0^Tg^*(t)f(t)\,dt\right)^n
	=\exp\left[\int_0^T(g^*(t)f(t)-1)\,dt\right],
$$
hence $\langle e(f),e(f)\rangle=e^{\|f\|^2_2-T}<\infty$ for any $f\in
L^\infty([0,T])$. We define $\mathsf{D}$, the exponential domain, as the
linear span of all $e(f)$, $f\in L^\infty([0,T])$, and we note that
$\mathsf{D}$ is dense in $\mathsf{F}$.  The exponential vectors have the
important property that they factorize over the continuous tensor product
structure (\ref{eq:guichardet}): indeed, it is evident from
(\ref{eq:exponential}) that $e(f)=e(f_{s]})\otimes e(f_{[s,t]})\otimes
e(f_{[t})$ where $f_{t]}$ is the restriction of $f$ to $[0,t]$, etc.

We are now ready to define a Poisson process.  Let us first define it as a
random variable on $\Omega$; we simply write $N_t(\tau)=|\tau\cap [0,t]|$, 
where $|\tau|$ denotes the number of elements in the set $\tau\in\Omega$.  
The random variable $N_t$ counts the number of jumps up to time $t$, and 
hence $\{N_t\}$ is by construction a Poisson process with unit rate under 
the measure ${\bf P}$.  We now turn this into an operator process by 
pointwise multiplication:
\begin{equation}
	(\Lambda_t\psi)(\tau)=N_t(\tau)\,\psi(\tau)=
	|\tau\cap [0,t]|\,\psi(\tau),
	\qquad \psi\in\mathsf{F},~\tau\in\Omega,~t\in[0,T].
\end{equation}
$\{\Lambda_t\}$ is called the {\it gauge process}; it is not difficult to 
see that though $\Lambda_t$ is an unbounded operator\footnote{
	As can be verified by explicit computation, the domain of 
	$\Lambda_t$ contains at least $\mathsf{D}$, the exponential 
	domain.  The reader may ask himself why we have only defined 
	exponential vectors $e(f)$ for $f\in L^\infty([0,T])$ rather
	than $f\in L^2([0,T])$: this is because the latter may not be in
	the domain of $\Lambda_t$.  Our domain $\mathsf{D}$ is sometimes
	called the {\it restricted} exponential domain in the literature.
}, it is affiliated to $\mathscr{W}_{t]}$ and hence the gauge process is 
adapted; in fact, the increments $N_t-N_s$ are even affiliated to 
$\mathscr{W}_{[s,t]}$.  Furthermore, $\Lambda_s$ and $\Lambda_t$ commute 
for all $s,t\in[0,T]$, and indeed ${\rm vN}(\Lambda_t,~t\in[0,T])=
L^\infty(\Omega,\mathcal{F},{\bf P})\subset\mathscr{W}$ is commutative.
Hence we could use the spectral theorem to map $\Lambda_t$ back to a 
classical stochastic process.  It is somewhat futile to diagonalize 
the operators using the spectral theorem, however, as we have already 
constructed them in diagonal form.

We have yet to introduce a state; a particularly interesting class of 
states are the {\it coherent states} $\mathbb{P}_f(X)=
\langle e(f),X\,e(f)\rangle\,e^{T-\|f\|_2^2}$.  Because of the continuous 
tensor product property, the coherent states split up as follows:
\begin{equation}\label{eq:indepincr}
	X=X_{s]}\otimes X_{[s,t]}\otimes X_{[t},\qquad
	\mathbb{P}_f(X)=\mathbb{P}_{f_{s]}}(X_{s]})~
	\mathbb{P}_{f_{[s,t]}}(X_{[s,t]})~
	\mathbb{P}_{f_{[t}}(X_{[t}).
\end{equation}
But as $N_t-N_s$ is affiliated to $\mathscr{W}_{[s,t]}$, it follows that 
under the state $\mathbb{P}_f$ the gauge process has independent 
increments.  Furthermore, if we denote by $P_{N_t-N_s}(B)$ the spectral 
measure of $N_t-N_s$, then we have
$$
	\mathbb{P}_f(P_{N_t-N_s}(B))=
	\mathbb{P}_{f_{[s,t]}}(\chi_B(|\tau\cap[s,t]|))=
	\sum_{n\in B}\frac{e^{-\int_s^t|f(r)|^2\,dr}}{n!}
		\left(\int_s^t|f(r)|^2\,dr\right)^n.
$$
Evidently, $\Lambda_t$ is an inhomogeneous Poisson process with rate
$|f(t)|^2$ under the state $\mathbb{P}_f$.  Note in particular that as
$e(1)(\tau)=1$, we have for any $X\in L^\infty(\Omega,\cF,{\bf P})$ the
relation $\mathbb{P}_1(X)=\langle 1,X\,1\rangle=E_{\bf P}(X)$; hence the
fact that under $\mathbb{P}_1$ the gauge process is a Poisson process with
unit rate is exactly what we expect from the definition of ${\bf P}$.  
Under $\mathbb{P}_0$, on the other hand, the gauge process doesn't
register any counts; $\mathbb{P}_0=\phi$ is called the {\it vacuum state},
and $e(0)=\Phi$ is called the {\it vacuum vector}.

\subsection{Weyl operators and Wiener processes}

We have now exhausted the diagonal observables affiliated to the space
$(L^\infty(\Omega,\mathcal{F},{\bf P}),\mathbb{P}_f)$: every such 
observable is some functional of the Poisson process $\Lambda_t$ with rate 
$|f|^2$.  Let us thus explore whether we can find interesting observables 
affiliated to $\mathscr{W}$ that do not commute with $\Lambda_t$.  To this 
end, we follow again essentially Example \ref{ex pos}.  Given $f,g\in 
L^\infty([0,T])$ we look for a unitary operator $W(f)$ that implements the 
translation group $W(f)e(g)\propto e(f+g)$.  A calculation identical to 
the one in Example \ref{ex pos} shows that we should define
\begin{equation}\label{eq Weyl}
  W(f)e(g) = 
	e^{-\int_0^T\left(f^*(t)g(t)
		+\frac{1}{2}f^*(t)f(t)\right)\,dt 
	} e(f+g)=e^{-\langle f,g\rangle_2-\|f\|^2_2/2}\,e(f+g).
\end{equation} 
The unitary operator $W(f)$ is called a {\it Weyl operator}, and provides
a projective unitary representation in the sense that $W(f)W(g)=
W(f+g)\,e^{i\,{\rm Im}\langle g,f\rangle_2}$.  Note that it is
sufficient to define the action of $W(f)$ only on exponential vectors; we
can extend to $\mathsf{D}$ by linearity, and as $\mathsf{D}$ is dense and
$W(f)$ is bounded the Weyl operators are uniquely extended to all of
$\qF$.  An important property, which follows immediately from the 
definition of $W(f)$ and the continuous tensor product property, is that
\begin{equation}\label{eq:weylsplit}
	W(f)e(g)=
		W(f_{s]})e(g_{s]})\otimes 
		W(f_{[s,t]})e(g_{[s,t]})\otimes 
		W(f_{[t})e(g_{[t}).
\end{equation}
In particular, we see that $W(f\chi_{[0,t]})$ is an adapted 
operator process.

Now fix $f\in L^\infty([0,T])$ and consider the unitary group
$\{W(tf)\}_{t\in\R}$; this group is in fact continuous \cite{Par92}, and
hence by Stone's Theorem \ref{thm stone} there exists a self-adjoint
$B(f)$ such that $W(kf)=e^{ikB(f)}$.  The operators $B(f)$, $f\in
L^\infty([0,T])$, are called \emph{field operators}.  Finding the
distribution of the observable $B(f)$ is straightforward, as the
characteristic function of $B(f)$ (under the coherent state
$\mathbb{P}_g$) is given by
$$
	b_f(k)=\mathbb{P}_g(W_{kf})
	=\langle e(g),e(g+kf)\rangle e^{T-\|g\|^2_2
		-k\langle f,g\rangle_2-k^2\|f\|^2_2/2}
	=e^{2ik\,{\rm Im}\langle g,f\rangle_2-k^2\|f\|^2_2/2}.
$$
Hence $B(f)$ is a Gaussian random variable with mean $2\,{\rm 
Im}\langle g,f\rangle_2$ and variance $\|f\|^2_2$.   In the vacuum, 
i.e.\ $g=0$, the mean vanishes; for simplicity, we will restrict ourselves 
to the vacuum case in the following.

Consider the operator process $\{B_t^\varphi=
B(e^{i\varphi}\chi_{[0,t]}):t\in[0,T]\}$ for some fixed, real function 
$\varphi\in L^\infty([0,T])$.  $B_t^\varphi$ is adapted, as we have
already established that $W(f\chi_{[0,t]})$ is adapted for any $f$; 
moreover, $B(e^{i\varphi}\chi_{[s,t]})=B_t^\varphi-B_s^\varphi$ is 
affiliated to $\mathscr{W}_{[s,t]}$ due to Eq.\ (\ref{eq:weylsplit}).  
This immediately tells us two important things: first, $B_t^\varphi$ and 
$B_s^\varphi$ commute for all $s,t\in[0,T]$; indeed, 
$B_t^\varphi-B_s^\varphi$ must commute with $B_s^\varphi-B_0^\varphi$, and
commutativity follows from $B_0^\varphi=I$.  This means that ${\rm
vN}(B_t^\varphi,~t\in[0,T])$ is a commutative algebra and hence we can
represent $B_t^\varphi$ for every $t$ as a classical random variable on 
the same probability space $(\Omega^\varphi,\mathcal{F}^\varphi,{\bf 
P}^\varphi)$; in particular, $\iota(B_t^\varphi)$ is a classical 
stochastic process.  Second, Eq.\ (\ref{eq:indepincr}) implies that the 
process $B_t^\varphi$ has independent increments.  But we have established  
$B_t^\varphi-B_s^\varphi$ is (in the vacuum) a mean zero Gaussian random 
variable with variance $t-s$, and as $B_t^\varphi$ has independent 
increments we have established that $\iota(B_t^\varphi)$ is precisely a 
Wiener process on $(\Omega^\varphi,\mathcal{F}^\varphi,{\bf P}^\varphi)$.

Let us introduce the following notation.  Define $Q_t=B(i\chi_{[0,t]})$,
$P_t=B(-\chi_{[0,t]})$, and $A_t=(Q_t+iP_t)/2$.  Note that $Q_t$ and $P_t$ 
are self-adjoint by Stone's theorem, whereas $A_t$ has the adjoint
$A_t^*=(Q_t-iP_t)/2$.  We now compute
$$
	B(f)e(g)=
	\left.\frac{1}{i}\frac{d}{dk}W(kf)e(g)\right|_{k=0}
	=i\langle f,g\rangle_2\,e(g)-i
	\left.\frac{d}{dk}e(g+kf)\right|_{k=0}.
$$
Evidently $A_te(g)=\langle\chi_{[0,t]},g\rangle_2\,e(g)=
\int_0^tg(s)ds\,e(g)$.  But then we can write
$$
	(A_te(g))(\tau)=\int_0^tg(s)\,ds\prod_{r\in\tau}g(r)=
	\int_0^tg(s)\prod_{r\in\tau}g(r)\,ds=
	\int_0^te(g)(\tau\cup\{s\})\,ds.
$$
In particular, this formula extends to any $\psi\in\qF$ for which the 
integral on the righthand side (with $e(g)$ replaced by $\psi$) defines a 
normalizable vector.  $A_t$ is called the Fock space {\it annihilation 
operator}, as it generalizes the corresponding notion introduced in 
Example \ref{ex pos}.  The reader should verify that its adjoint can be 
expressed as
$$
	(A^*_t\psi)(\tau)=\sum_{s\in\tau\cap[0,t]}\psi(\tau
		\backslash\{s\})
$$
on a sufficiently large domain.  Not surprisingly, $A^*_t$ is called the 
{\it creation operator}.  It is conventional in quantum stochastic 
calculus to use $A_t$ and its adjoint rather than $Q_t$ and $P_t$; we 
shall conform to this standard.

In summary, we have constructed a quantum probability space $({\mathscr
W},\phi)$ that admits an entire family (indexed by $\varphi$) of Wiener
processes.  Note however, that these processes do not necessarily commute
for different $\varphi$; in fact, it is not difficult to establish that
$[B(f),B(g)]\psi=2i\,{\rm Im}\langle f,g\rangle_2\,\psi$ on a suitably
large domain (e.g.\ $\psi\in\mathsf{D}$).  Therefore, even though every
$B_t^\varphi$ defines a Wiener process, these cannot be represented on the
same classical probability space for different $\varphi_{1,2}$ unless 
${\rm Im}(e^{i(\varphi_1-\varphi_2)})=0$.

We have also defined a Poisson process $\Lambda_t$, but unfortunately it
vanishes in the vacuum.  Consider, however, the process
$\Lambda_t(f)=W(f)^*\Lambda_tW(f)$; for any Borel function $b$ we can
write $\phi(b(\Lambda_{t_1}(f),\ldots,\Lambda_{t_n}(f)))=
\mathbb{P}_f(b(\Lambda_{t_1},\ldots,\Lambda_{t_n}))$.  Evidently
$\Lambda_t(f)$ has the same statistics in the vacuum as does $\Lambda_t$
under the coherent state $\mathbb{P}_f$.  This shows that we can define
even a whole family of Poisson processes in the vacuum.  We do not lose
much by restricting ourselves to the vacuum as an underlying state (as we
will do in the remainder of the article), as we can always transform to a
coherent state by ``sandwiching'' with Weyl operators.  Note that like the
family $B_t^\varphi$, the processes $\Lambda_t(f)$ do not commute amongst
each other.  We see that the quantum probability space
$(\mathscr{W},\phi)$ gives rise to a rich family of incompatible
stochastic processes.

\subsection{Quantum stochastic calculus}\label{sec qsc}

Now that we have obtained Wiener and Poisson processes, we can try to
develop stochastic integrals with respect to these processes and an
associated stochastic calculus.  Note that if we were only interested in,
e.g., integrating with respect to $Q_t$ an adapted process which commutes
with $Q_t$, then we could simply use the classical It\^o integral
definition through the spectral theorem.  This will not suffice for our
purposes, however, as we will want to consider stochastic differential
equations that are driven simultaneously by the noncommuting noises $Q_t$
and $P_t$ (and even $\Lambda_t$).  Moreover, we would like to have an
It\^o rule that tells us how to multiply stochastic integrals with
respect to $Q_t$ and $P_t$.

Our motivation for developing generalized quantum stochastic calculus is
that this allows us to rigorously define and manipulate Schr{\"o}dinger
equations, as in Eq.\ (\ref{schrodinger}), with a white-noise Hamiltonian
formally defined by $H(t)=H_0+H_1\,\dot Q_t+ H_2\,\dot P_t$.  In
\S\ref{sec filtering} we will see that such models emerge naturally in
applications.  In this section we sketch the development of quantum
stochastic calculus as it was introduced in a seminal paper by Hudson and
Parthasarathy \cite{HuP84}.  For a full development of this calculus we 
refer to \cite{HuP84,Hud03,Par92}. The Hudson-Parthasarathy approach has 
some technical issues, not surprisingly involving the unboundedness of 
operators, the full extent of which is still being explored.  Though we 
cannot go into detail here, we will attempt to sketch some of the issues 
and give references to recent literature.

We work in the following setting.  We wish to integrate processes against
the three noises $A_t$, $A_t^*$ and $\Lambda_t$ (the {\it fundamental
noises}), i.e.\ we want to define $\int_0^tL_s\,dM_s$ where $M_t$ is one
of the fundamental noises.  The noises are defined on the quantum
probability space $(\mathscr{W},\phi)$, but we will want to couple these
noises to an external quantum system, the {\it initial system}\footnote{
	This name has the following origin.  Recall from \S\ref{sec BG} 
	that observables $X$ evolve in time as $X_t=U_t^*XU_t$ (we will 
	define a unitary evolution $U_t$ in \S\ref{sec filtering}).
	We would like to think of $X\otimes I\in\mathscr{B}\otimes
	\mathscr{W}$ as describing the external system; however,
	$U_t^*(X\otimes I)U_t$ will not be of the form $Y\otimes I$ except 
	at $t=0$. Hence the initial system observable $X\otimes I$ 
	describes the external system at the initial time $t=0$.
}, with which they interact.  To this end, let us introduce the initial
Hilbert space $\qh$, $\mathscr{B}=\mathscr{B}(\qh)$ and the associated
initial quantum probability space $(\mathscr{B},\rho)$.  We will choose
our integrands $L_t$ to be adapted processes on
$(\mathscr{B}\otimes\mathscr{W}, \rho\otimes\phi)$, i.e.\ each $L_t$ is
affiliated to $\mathscr{B}\otimes\mathscr{W}_{t]}$ and acts as $I$ on
$\mathscr{W}_{[t}$.

As usual, we begin with simple processes.  Given $s<t$, recall that for
the fundamental processes $M_t-M_s$ is affiliated to
$\mathscr{W}_{[s,t]}$, whereas for adapted processes $L_s$ is affiliated
to $\mathscr{W}_{s]}$; hence we can naturally write
$L_s(M_t-M_s)=L_s\otimes(M_t-M_s)$. In particular the increment $M_t-M_s$
commutes with $L_s$, and we have no problems with operator multiplication
of these unbounded operators.  Let $\{t_i:i=0,\ldots,n,~t_i<t_{i+1}\}$ be
a sequence of times with $t_0=0$ and $t_n=T$.  By definition, we set 
$$
	L_t=\sum_{i=0}^{n-1}L_{t_i}\chi_{[t_i,t_{i+1})}(t)
	\quad
	\Longrightarrow
	\quad
	\int_0^t L_s\,dM_s=\sum_{i=0}^{n-1}
	L_{t_i}\otimes (M_{t_{i+1}\wedge t}-M_{t_i\wedge t}).
$$
This definition makes sense as long as the operators $L_t$ and $M_t$ have
a sufficiently large common dense domain that the sum is well defined.  
To enforce this, we will require that the domain of every $L_t$ contains 
at least the exponential domain $\mathsf{D}$.

Now comes the hard part in any integration theory: given a quadruple 
of suitably restricted adapted processes $(E,F,G,H)$, such that these 
admit simple approximations $(E^{n},F^{n},G^{n},H^{n})$, we 
wish to define the integral
\begin{equation}\label{eq:stochintgen}
	I_t=\int_0^t(E_t\,d\Lambda_t
		+F_t\,dA_t+G_t\,dA_t^*+H_t\,dt)
\end{equation}
as a limit, in some sense, of the corresponding integrals $I_t^{n}$ over 
the simple processes.  Recall that in the classical theory, the 
It\^o isometry allows us to define the stochastic integral as a 
mean-square limit of simple processes, and a little more work shows that 
every square-integrable process admits a mean-square approximation by 
simple processes.  Things are not quite so ``simple'' in the 
noncommutative case, however.

To see what goes wrong, consider for simplicity the case $\qh={\bf C}$ so
that we can forget about the initial state $\rho$.  We already encountered
the noncommutative $L^2$ (semi)norm $\|X\|_\phi^2=\phi(X^*X)$ when we
discussed conditional expectations.  We are thus looking for a suitable
unbounded operator $I_t$ such that we have mean-square convergence,
$\|I_t-I_t^{n}\|_\phi^2= \langle
(I_t-I_t^{n})\Phi,(I_t-I_t^{n})\Phi\rangle\to 0$ as $n\to\infty$.  
But this is a very ill-defined problem, as it only depends on the action
of $I_t$ on the vacuum vector $\Phi$; in particular, what do we choose as
the domain of $I_t$, and how do we define $I_t$ on vectors orthogonal to
$\Phi$?  There could be a large number of inequivalent ways of doing this,
giving rise to limiting operators with very different properties\footnote{
	This was not a problem for the definition of conditional 
	expectations;  as all versions of the conditional expectation are
	affiliated to a single commutative algebra, they are
	a.s.\ equivalent by the spectral theorem.  On the other hand, 
	various ``versions'' of $I_t$ that satisfy 
	$\|I_t-I_t^{n}\|_\phi\to 0$ need not even commute, and such 
	operators are fundamentally inequivalent.
}.

The solution of Hudson and Parthasarathy works as follows.  First of all,
we fix the domain of $I_t$ at the outset: every stochastic integral will
have $\qh\otimes\mathsf{D}$ as its domain (one could choose a dense domain
in $\qh$ as well; we will not worry about this).  To specify $I_t$ as a
limit of simple integrals $I_t^n$, we choose $I_t$ as the unique operator
on $\qh\otimes\mathsf{D}$ such that $\langle (I_t-I_t^n)\,v\otimes\psi,
(I_t-I_t^n)\,v\otimes\psi\rangle\to 0$ for every $\psi\in\mathsf{D}$,
$v\in\qh$ (it is sufficient to verify this for $\psi=e(f)$, $f\in
L^\infty([0,T])$).  In essence this is like a mean-square limit, but
simultaneously for every coherent state.  A suitable estimate replaces the
It\^o isometry \cite[Corollary 1]{HuP84} and shows that this limit exists
as long as $\int_0^T\|(E_s-E_s^n)\,v\otimes\psi\|^2ds\to 0$ as
$n\to\infty$ for every $\psi\in\mathsf{D}$, $v\in\qh$ (and similarly for
$F,G,H$), independent of the approximation.  Finally, \cite[Proposition
3.2]{HuP84} shows that every square-integrable process, i.e.\
$\int_0^T\|E_s\,v\otimes\psi\|^2ds<\infty$ for all $\psi\in\mathsf{D}$,
$v\in\qh$, admits a suitable approximation by simple processes.  We thus
arrive at the following.

\begin{definition}[Quantum It\^o integral]
	An operator process $\{X_t\}$ is stochastically integrable if 
	it is adapted and square-integrable.
	Given a quadruple $(E,F,G,H)$ of such processes, the stochastic
	integral (\ref{eq:stochintgen}) is uniquely defined as the limit
	of simple approximations on the domain $\qh\otimes\mathsf{D}$.
\end{definition}

A property that we will exploit in future is $\Lambda_t\Phi=A_t\Phi=0$.  
It is immediate from the definition that stochastic integrals with respect
to $A_t$ and $\Lambda_t$ acting on $\Phi$ vanish.  Hence the vacuum
expectations of stochastic integrals with respect to $A_t$ and $\Lambda_t$
vanish as well.  Furthermore, as $\langle\Omega,A^*_t\Omega
\rangle=\langle A_t\Omega,\Omega\rangle=0$, we see that at least for
simple processes (and indeed this holds for any integrand) the vacuum
expectation of stochastic integrals with respect to $A^*_t$ vanish. Note,
however, that $A^*_t\Phi\neq 0$.

Our next task is to develop a stochastic calculus; the integrals defined
above are not of much use, unless we have an It\^o product rule with which
they can be manipulated.  Once again we run into unpleasant problems. If
$I_t$ and $J_t$ are integrals of the form (\ref{eq:stochintgen}), there is
no reason to expect that their product $I_tJ_t$ is a well-defined operator
on the domain $\qh\otimes\mathsf{D}$.  The idea of Hudson and
Parthasarathy is inspired by the identity $\langle\psi',X^*Y\psi\rangle
=\langle X\psi',Y\psi\rangle$ for bounded operators; rather than finding
an expression for $I_tJ_t$, they calculate $\langle I_t\,v'\otimes\psi',
J_t\,v\otimes\psi\rangle$ for every $v\in\qh$, $\psi\in\mathsf{D}$, which
is always well defined.  One finds explicitly a lengthy expression
\cite[Theorems 4.3--4.4]{HuP84}, which is essentially the quantum It\^o
rule expressed in terms of $\qh\otimes\mathsf{D}$-matrix elements.

In practice, however, we are mostly interested in calculating actual 
operator products $I_tJ_t$.  We will need the concept of an {\it adjoint 
pair}; two operators $X$ and $X^\dag$ are said to be an adjoint pair if 
$\langle v'\otimes\psi',X\,v\otimes\psi\rangle=\langle 
X^\dag\,v'\otimes\psi',v\otimes\psi\rangle$ for every $v\in\qh$, 
$\psi\in\mathsf{D}$.  It is not difficult to verify that if $(E,F,G,H)$ 
and $(E^\dag,F^\dag,G^\dag,H^\dag)$ are adjoint pairs, then $I_t$ and 
$I_t^\dag$ form an adjoint pair, where
\begin{equation}\label{eq:stochintadjoint}
	I_t^\dag=\int_0^t(E_t^\dag\,d\Lambda_t
		+F_t^\dag\,dA_t^*+G_t^\dag\,dA_t+H_t^\dag\,dt).
\end{equation}
In essence, the adjoint $\dag$ replaces the Hilbert space adjoint $*$ on
the domain $\qh\otimes\mathsf{D}$.  Now suppose that we can verify
explicitly that the product $I_tJ_t$ is well defined; then we can read off
an expression for $I_tJ_t$ from the matrix elements $\langle
I_t^\dag\,v'\otimes\psi', J_t\,v\otimes\psi\rangle$.  This gives the
following explicit form of the quantum It\^o rule.

\begin{theorem}[Quantum It\^o rule {\cite[Proposition 25.26]{Par92}}]
\label{Itorule}
Let $(F,G,H,I)$, $(B,C,D,E)$ and $(B^\dag,C^\dag,D^\dag,E^\dag)$ be 
quadruples of stochastically integrable processes such that the latter two 
quadruples are adjoint pairs.   Define the stochastic integrals
\begin{equation*}\begin{split}
	& dX_t = B_t\,d\Lambda_t + C_t\,dA_t + D_t\,dA_t^* + E_t\,dt, \\
	& dY_t = F_t\,d\Lambda_t + G_t\,dA_t + H_t\,dA_t^* + I_t\,dt,
\end{split}\end{equation*}
and suppose that we have verified that the product $X_tY_t$ is well 
defined and that $X_tF_t,\ldots,X_tI_t$, $B_tY_t,\ldots,E_tY_t$, and
$B_tF_t,B_tG_t,\ldots,E_tI_t$ are well defined and stochastically 
integrable.  Then the process $X_tY_t$ satisfies the relation
$$
	d(X_tY_t) = X_t\,dY_t + (dX_t)\,Y_t + dX_t\,dY_t,
$$
where $X_t\,dY_t=X_tF_t\,d\Lambda_t + X_tG_t\,dA_t + X_tH_t\,dA_t^* + 
X_tI_t\,dt$, $(dX_t)\,Y_t=B_tY_t\,d\Lambda_t + C_tY_t\,dA_t + 
D_tY_t\,dA_t^* + E_tY_t\,dt$, and $dX_t\,dY_t = B_tF_t\,d\Lambda_t + 
C_tF_t\,dA_t + B_tH_t\,dA^*_t + C_tH_t\,dt$ is evaluated according to 
the quantum It\^o table
\begin{center}
\vskip.2cm
{\begin{tabular}{l|llll}
	dX$~\,\backslash~$dY
		     & $dA_t$ & $d\Lambda_t$ 	& $dA^*_t$ 	 & $dt$\\
\hline 
	$dA_t$ 	     & $0$    & $dA_t$ 		& $dt$ 		 & $0$ \\
	$d\Lambda_t$ & $0$    & $d\Lambda_t$ 	& $dA^*_t$ 	 & $0$ \\
	$dA^*_t$     & $0$    & $0$  		& $0$ 		 & $0$ \\
	$dt$ 	     & $0$    & $0$ 		& $0$ 		 & $0$
\end{tabular}}
\vskip.2cm
\end{center}
In particular, the theorem holds if $B_t,C_t,D_t,E_t$ and $X_t$ are
bounded processes {\rm \cite{HuP84}}, in which case the adjoints $B^\dag$ 
etc.\ are simply taken to be the Hilbert space adjoints $B^*$ etc., and 
$X_t$ extends uniquely to a bounded operator in $\mathscr{W}_{t]}$. 
\end{theorem}

\begin{remark}
The choice to restrict attention to a fixed domain $\qh\otimes\mathsf{D}$ 
allows Hudson-Parthasarathy to develop a viable quantum stochastic 
calculus.  This choice, however, has quite a few drawbacks; we highlight 
one of the problems.  Suppose $X$ is self-adjoint; implicit in this 
statement is that the domains of $X$ and $X^*$ coincide.  It can 
happen that if we restrict the domain of $X$, then the restricted operator 
admits many inequivalent self-adjoint extensions; see 
\cite[pages 257--259]{ReS80} for an example.  Hence the restriction to a 
fixed domain can become a real, physical problem, that prevents us from 
uniquely interpreting unbounded operators on $\qh\otimes\mathsf{D}$ as 
observables.

Such problems have prompted the development of alternative approaches to
quantum stochastic integration, and the topic is still under active
investigation.  In a significant recent achievement Attal and Lindsay,
building on several earlier approaches (see e.g.\ \cite{Mey93,Bia95} 
and the references therein), develop a theory in which
the integrals achieve their maximal domains \cite{AtL04}.  Unfortunately,
the theory is very technical and a little daunting for every-day use.  A
different approach that even preceeds Hudson and Parthasarathy is that of
Barnett, Streater and Wilde \cite{BSW83}.  Their theory is attractive as 
it is completely algebraic in nature (the Hilbert space and its domains do 
not play a fundamental role), but lacks a satisfactory It\^o rule.

Despite these issues, the Hudson-Parthasarathy approach works quite well.  
In practice one usually works with a ``noisy Schr{\"o}dinger equation''
Eq.\ (\ref{eq Ut}), the solution of which is unitary and thus bounded.  
As long as the integrals and integrands are bounded, they are uniquely
defined by their specification on a dense domain.  In this article, in
keeping with our attitude towards unbounded operators, we will not
worry about such issues and assume that we can apply the quantum It\^o
rules. 
\qed 
\end{remark}

\begin{example}\label{ex qsc}
In \S\ref{sec filtering} we will encounter quantum stochastic differential 
equations (QSDE), the treatment of which proceeds along the same lines as 
the classical theory.  We claim that the Weyl operator $W(f_{t]})$ is the 
solution of the QSDE
\begin{equation}\label{eq dif Weyl}
	dW(f_{t]}) = \Big\{
		f(t)\,dA^*_t-f^*(t)\,dA_t - \frac{1}{2}|f(t)|^2\,dt
	\Big\}W(f_{t]}).
\end{equation}
In particular, one can verify the Weyl relation $W(f)W(g)=
W(f+g)\,e^{i\,{\rm Im}\langle g,f\rangle_2}$ directly using the 
quantum It\^o rule.  From Eq.\ (\ref{eq dif Weyl}) and $W(kf)=e^{ikB(f)}$ 
we obtain
$$
	B(f)=\int_0^T(if(t)^*\,dA_t-if(t)\,dA_t^*).
$$
Hence $dB_t^\varphi=ie^{-i\varphi(t)}\,dA_t-ie^{i\varphi(t)}\,dA_t^*$, and 
the quantum It\^o rules reduce to the classical It\^o rule
$(dB_t^\varphi)^2=dt$.  Finally, recall that we defined Poisson processes 
$\Lambda_t(f)=W(f)^*\Lambda_t W(f)=W(f_{t]})^*\Lambda_tW(f_{t]})$ (the 
latter equality is due to $W(f)=W(f_{t]})\otimes W(f_{[t})$ and the fact 
that $W(f_{[t})\in\mathscr{W}_{[t}$ is unitary and commutes with the 
adapted process $\Lambda_t$). Using the quantum It\^o rule we obtain the 
explicit representation
\begin{equation}\label{eq:poissonexplicit}
	d\Lambda_t(f) = d\Lambda_t + f^*(t)\,dA_t + f(t)\,dA_t^* + 
	|f(t)|^2\,dt, 
\end{equation}
for which the quantum It\^o rules reduce to the classical product rule 
$(d\Lambda_t(f))^2=d\Lambda_t(f)$ for a Poisson process.
\qed
\end{example}


\section{The filtering problem in quantum optics}
\label{sec filtering}

Many realistic physical scenarios are very well described by quantum
stochastic differential equations driven by the processes $A_t$,
$A^*_t$ and $\Lambda_t$ discussed in the previous section.  Of
course, as in the classical theory, white noise systems are only an
idealization of physical interactions; a Markov limit of wide-band noise
in the spirit of Wong and Zakai \cite{Gou05} gives stochastic models in
the It\^o form.  For a large class of quantum systems, particularly those
arising in the field of quantum optics, such approximations are extremely
good and describe laboratory experiments essentially to experimental
precision.  Though a detailed discussion of the physics involved in the
modelling of such systems is beyond the scope of this article, we here 
very briefly describe the physical origin of the equations that are widely 
used in the physics community \cite{GC84}, describe the measurements that 
are made, and set up the quantum filtering problem to be solved.

\subsection{The quantum optics model}
\label{sec filtering model}

The basic model of quantum optics consists of some fixed physical system,
e.g.\ a collection of atoms, in interaction with the electromagnetic
field. The atomic observables are self-adjoint operators on a Hilbert
space $\qh$.  The description of the electromagnetic field and its
interaction with the atoms follows from basic physical arguments (see the
excellent monograph \cite{CT89} for a thorough treatment of this theory,
known as {\it quantum electrodynamics}).  It turns out that the free
electromagnetic field, i.e.\ an optical field in empty space, is described
by a stationary Gaussian (noncommutative) wide band noise $\tilde a(t,{\bf
r})$ that propagates through space at the speed of light $c$; i.e.\ if we
restrict ourselves to a single spatial dimension, $\tilde
a(t+\tau,z)=\tilde a(t,z-c\tau)$.  If we now place the atoms at the origin
$z=0$, then the quantum dynamics is given by a Schr\"odinger equation of
the form
\begin{equation}\label{eq:colorednoise}
	\frac{d}{dt}\tilde U(t)=
	\left[-iH+L\,\tilde a^*(t,0)-L^*\,\tilde a(t,0)
	\right]\tilde U(t),\qquad U(0)=I,
\end{equation}
where $L\in\mathscr{B}$ is an atomic (dipole) operator and 
$H\in\mathscr{B}$ is an atomic Hamiltonian, $H$ being self-adjoint.  This 
equation, which follows directly from the physical model, has wide-band 
right hand side. Note that we have set $\hbar=1$ for convenience, a 
convention ubiquitous in physics (the only consequence is a change of 
units).

We now want to approximate the wide-band noise by white noise.  This can 
be done in a rigorous way \cite{AFLu90,AGLu95,Gou05}, but we will not 
detail the procedure here (a brief sketch can be found in \cite{HSM05a}).
Suffice it to say that one arrives at the following {\em quantum 
stochastic differential equation} (QSDE)
\begin{equation}\label{eq Ut}
	dU_t = \Big\{L\,dA_t^* -L^*\, dA_t - \frac{1}{2}L^*L\,dt - 
		iH\,dt\Big\}U_t, \qquad U_0 = I,
\end{equation}
which is driven by the non-commuting white noise processes $A_t$ and
$A^*_t$.  Note that this is almost precisely of the same form as Eq.\ 
(\ref{eq:colorednoise}), except that we have added the It\^o correction 
term $-\tfrac{1}{2}L^*LU_t\,dt$.  A Picard iteration argument 
\cite{HuP84,Par92} ensures existence and uniqueness of the solution. 
The adjoint $U^*_t$ satisfies
\begin{equation*}
	dU^*_t = U^*_t\Big\{L^*\, dA_t -L\,dA^*_t - 
	\frac{1}{2}L^*L\,dt + iH\,dt\Big\}, \qquad U^*_0 = I.
\end{equation*}
Using the quantum It\^o rule we can calculate $d(U^*_tU_t)= 
d(U_tU_t^*) = 0$, i.e.\ the solution $U_t$ is unitary for 
all $t$ (as the solution of a Schr\"odinger equation should be).

Henceforth we will take Eq.\ (\ref{eq Ut}) as our physical model.  $U_t$ 
defines the time evolution or {\it flow} $j_t:X\mapsto U_t^*(X\otimes 
I)U_t$ of every atomic observable $X\in\mathscr{B}$ (recall the time 
evolution in \S\ref{sec BG QM}); i.e., an observation of $X \in 
\mathscr{B}$ at time $t$ is described by the observable $X_t=j_t(X)$.  
Using the It\^o rules, we find an explicit dynamical equation
\begin{equation}
	dj_t(X)=j_t(\mathcal{L}_{L,H}(X))\,dt+
	j_t([L^*,X])\,dA_t+j_t([X,L])\,dA_t^*, \ X \in \mathscr{B}, 
		\label{jt-qsde}
\end{equation}
where the so-called Lindblad generator \cite{Lin76} is given by
\begin{equation*}
	\mathcal{L}_{L,H}(X) = i[H,X] + L^*XL - 
	  \frac{1}{2}(L^*LX+XL^*L), \qquad X \in {\mathscr B}.
\end{equation*}   
In quantum probability, this object plays the same role as the 
infinitesimal generator of a Markov diffusion in classical probability 
theory.

\begin{remark}
Though it is unusual, one could use a very similar notation in classical
stochastic models.  Suppose some system is described by an underlying
configuration $x_t$ that obeys $dx_t=b(x_t)\,dt+ \sigma(x_t)\,dW_t$.  
Then the ``observables'' in the theory, i.e.\ things we could try to
measure, are functions $f$ of the configuration of the system.  The
observable $f$ at time $t$ is described by the random variable
$j_t(f)=f(x_t)$.  Using the classical It\^o rules, we get
$dj_t(f)=j_t(\mathcal{L}f)\,dt+j_t(\Sigma f)\,dW_t$ where
$\mathcal{L}f(x)=\sum_ib^i(x)\partial_if(x)+\tfrac{1}{2}
\sum_{ij}\sigma^i(x)\sigma^j(x)\partial_i\partial_jf(x)$ is the generator
of the Markov diffusion $x_t$, and $\Sigma f(x)=
\sum_i\sigma^i(x)\partial_if(x)$.  This expression is the classical analog
of (\ref{jt-qsde}); the sample paths $x_t$ do not have a quantum
counterpart, however.  \qed
\end{remark}

\subsection{Measurements}
\label{sec filtering measurements}

\begin{figure}
\centering
\includegraphics[width=0.9\textwidth]{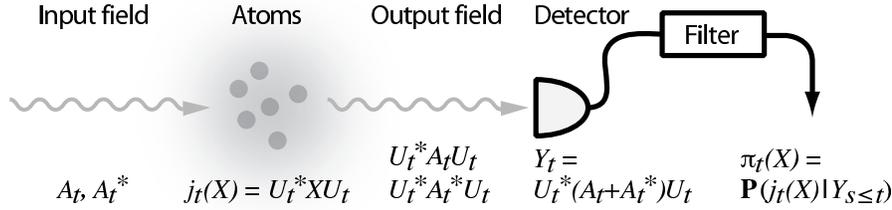}
\vskip-1cm
\caption{Cartoon of the quantum filtering setup in quantum optics.  An 
optical field, described by the field operators $A_t$, $A_t^*$, interacts 
with a system, e.g.\ a cloud of atoms.  After the atom-field interaction 
the field operators, as well as system operators $X$, are rotated by the 
unitary $U_t$.  The field is detected, giving rise to the observation 
$Y_t$.  Finally, the quantum filter (implemented on a classical signal 
processor) estimates atomic observables based on the field observations.
\label{fig:model}}
\end{figure}

Having described the system and its interaction with the field, let us now
turn to the observations that we can perform.  Unlike in classical models,
where one observes the system directly (with the addition of some
corrupting noise), in quantum models an observation is generally performed
in the field.  From the system's perspective, the interaction with the
field looks like an (albeit noncommutative) noisy driving force.  
Similarly, however, the field is perturbed by its interaction with the
atoms, and carries off information as it propagates away after the
interaction.  By performing a measurement in the field, then, we can
attempt to perform statistical inference of the atomic observables.  The
entire setup is depicted in Fig.\ \ref{fig:model}.

To calculate the perturbation of the field by the atoms we once again
calculate $U_t^* YU_t$, where now, however, $Y$ is a field observable. The
field observable of interest depends on the type of measurement we choose
to perform.  Without entering into the details, we mention two types of 
measurement that are extremely common in quantum optics:  direct 
photodetection (photon counting), for which the observation at time 
$t$ is given by $Y_t^\Lambda=U_t^*\Lambda_tU_t$, and homodyne 
detection, for which $Y_t^W=U_t^*(A_t+A_t^*)U_t$ (more 
generally $Y_t^W=U_t^*(e^{-i\varphi}A_t+e^{i\varphi}A_t^*)U_t$).  We 
refer to \cite{Bar90,Bar03} for a detailed treatment of quantum optical 
measurements.  Using the It\^o rules we obtain 
\be
	dY_t^\Lambda=d\Lambda_t+j_t(L)\,dA_t^*
		+j_t(L^*)\,dA_t+j_t(L^* L)\,dt,
	\label{y lambda}
\ee
\be
	dY_t^W=j_t(L+L^*)\,dt+dA_t+dA_t^*.
	\label{y w}
\ee
Intuitively, it would appear that $Y_t^\Lambda$ is like a Poisson process 
whose intensity is controlled by $j_t(L^*L)$ (recall Example 
\ref{ex qsc}), whereas $Y_t^W$ looks like a noisy observation of 
$j_t(L+L^*)$.  One should be careful with this conclusion, however, 
as $j_t(L)$ need not commute with $A_t$ or $A_t^*$, nor with itself 
at different times. 

It is essential, however, that the observation process commutes with 
itself at different times, and is hence equivalent to a classical 
stochastic process through the spectral theorem.  An observation process 
that does not obey this property cannot be observed in a single realization 
of an experiment and is physically meaningless.  Let us show that the 
observations processes we have defined above do obey this property, 
which is called the {\it self-nondemolition} property.  Let $Z$ be any 
operator of the form ${I}\otimes Z_{s]}\otimes{I}$ on
$\qh\otimes\qF_{s]}\otimes\qF_{[s}$ and let $t\ge s$.  Then the It\^o 
rules give directly
$$
	U_t^* ZU_t=U_s^* ZU_s+\int_s^t
	U_\tau^*\mathcal{L}_{L,H}(Z)U_\tau\,d\tau
	+\int_s^t U_\tau^* [L^*,Z]U_\tau\,dA_\tau
	+\int_s^t U_\tau^* [Z,L]U_\tau\,dA_\tau^*.
$$
Now let $Z=A_s+A_s^*$ or $Z=\Lambda_s$.  In both cases
$\mathcal{L}_{L,H}(Z)=[Z,L]=0$ as $L$ and $H$ are system observables and
$Z$ is a field observable.  Hence $Y_s^W=U_t^*(A_s+A_s^*)U_t$ and
$Y_s^\Lambda=U_t^*\Lambda_s U_t$ for all $t\ge s$.  It is now easily
verified, using the unitarity of $U_t$ and the fact that $A_s+A_s^*$ and
$\Lambda_s$ are commutative processes, that
$[Y_t^W,Y_s^W]=[Y_t^\Lambda,Y_s^\Lambda]=0$ for all $t,s$.  We denote by
${\mathscr Y}^W_t$ and ${\mathscr Y}^\Lambda_t$ the commutative von
Neumann algebras generated by the observation processes $Y_s^W$ and
$Y_s^\Lambda$, $s \leq t$, respectively. Do note, however, that $Y_t^W$
and $Y_t^\Lambda$ do not commute with each other; in any experiment, we
can choose to perform only one of these measurements.  Once we have made 
this choice, however, we can use the spectral theorem to represent the 
observations $Y_t$ as a classical stochastic process $\iota(Y_t)$ on a 
probablity space.

\subsection{Statement of the filtering problem}
\label{sec filtering statement}

Moving on to the next step in our program, we now wish to use the
information gained from the measurement process to infer something about
the system.  To find a least mean square estimate of a system observable 
$X\in\mathscr{B}$ at time $t$, given the observations $Y_t$ up to 
this time, we must calculate the conditional expectation
\be
	\pi_t(X)=\mathbb{P}(j_t(X)|\mathscr{Y}_t)
	\label{jt-X-cdl-exp}
\ee
where $\mathscr{Y}_t={\rm vN}(Y_s:0\le s\le t)$.  The remainder of this 
article is devoted to finding a recursive equation for $\pi_t(X)$ (the 
{\it filtering equation}).  Recall, however, that the conditional 
expectation is only defined if $j_t(X)$ is in the commutant of 
$\mathscr{Y}_t$, the interpretation being that statistical inference of an 
observable is only physically meaningful if the conditional statistics 
could possibly be tested through a compatible experiment.  Through an 
entirely identical procedure to the one used to show the 
self-nondemolition property, we can show that $j_t(X)$ is in the commutant 
of $\mathscr{Y}_t$ for any $X\in\mathscr{B}$.  This is known as the {\it 
nondemolition property}, which can be written as
\be
	[ j_t(X), Y_s ] = 0 \ \ \forall \ s \leq t , \  X \in \mathscr{B} .
	\label{qnd}
\ee
We note that we have now obtained a system-theoretic model of our system 
and observations, defined on the quantum probability space
$(\mathscr{B}\otimes\mathscr{W},\mathbb{P}=\rho\otimes\phi)$ by
\begin{eqnarray}
	& dj_t(X)=j_t(\mathcal{L}_{L,H}(X))\,dt+
		j_t([L^*,X])\,dA_t+j_t([X,L])\,dA_t^* ,
	\label{X} \\
	& dY_t=j_t(L+L^*)\,dt+dA_t+dA_t^*
	\label{Y}
\end{eqnarray}
in the case of homodyne detection, or by Eq.\ (\ref{X}) and
\be
	dY_t=d\Lambda_t+j_t(L)\,dA_t^*
		+j_t(L^*)\,dA_t+j_t(L^* L)\,dt
	\label{YL}
\ee
in the case of counting observations.  These equations define a 
system-observation model in direct analogy to such models used throughout 
classical nonlinear filtering and stochastic control theory. 

\begin{remark}  \label{rmk measurement noise}
Unlike in a classical filtering scenario, we have not added any
independent corrupting noise to the observations.  Nonetheless, the
filtering problem does not reduce to a problem with complete observations
because the system is driven by noise that does not commute with the 
observations.  Hence the problem of partial observations is intrinsic to 
quantum measurement theory.  The quantum filtering problem considered here 
is the simplest possible one; one could add additional corrupting noise as 
in the classical case, have the system interact with multiple fields (some 
of which are observed, others unobserved), etc.  These are not essential 
complications, however, and filters for such models are obtained much in 
the same way. \qed
\end{remark}


\section{The reference probability method}
\label{sec KS formula}

The goal of this section is to derive the quantum filtering equation, a 
recursive equation for $\pi_t(X)$, using a method that is close to the 
classical reference probability method of Duncan \cite{TD68}, Mortensen 
\cite{RM66}, and Zakai \cite{Zak69}.  We consider first the homodyne 
detection case, then the photon counting case.  In \S\ref{sec mtg} we will 
rederive the filtering equation for the homodyne detection case using 
martingale methods; the chief advantage of the reference probability 
method is that it is somewhat simpler to apply.  The following approach is 
based on \cite{BH05}.

\subsection{Homodyne detection}

Let us briefly recall the classical reference probability procedure; for
an introduction see e.g.\ \cite{AEM95}.  In order to simplify
the filtering problem, one starts by introducing a new probability
measure, using a Girsanov transformation, under which the measurement
record is a Wiener process.  Then various (elementary) properties of the
conditional expectation allow the filtering problem to be expressed, and
solved, with respect to the new measure.  We now apply this logic to the
quantum filtering problem. Note that we have already applied the method in 
Example \ref{eg stern};  the following is essentially a continuous time 
version of that example.

We consider the homodyne detection setup given by Eqs.\ (\ref{X}) and 
$(\ref{Y})$.  We could try to find a new state under which $Y_t$ is a 
Wiener process; however, it will be more convenient to work not in terms 
of $Y_t$ but in terms of $Z_t=A_t+A_t^*$, as it is very easy to manipulate 
$Z_t$ using the methods of \S\ref{sec Fock}.  Thus before we really start 
filtering, let us transform the problem in terms of $Z_t$.  Introduce the 
state $\mathbb{Q}^t$ defined by
\be
	\mathbb{Q}^t (X) = \mathbb{P}( U^*_t X U_t ) ,
	\label{schrodinger-state}
\ee
with $U_t$ as in \S\ref{sec filtering}, and we fix from now on
$\mathbb{P}=\rho\otimes\phi$.
Now recall from Example \ref{eg stern} that $\mathbb{Q}(X)=
\mathbb{P}(U^*XU)$ implies $\mathbb{P}(U^*XU|U^*\mathscr{C}U)=
U^*\mathbb{Q}(X|\mathscr{C})U$ (this is easily checked using the 
definition of the conditional expectation).  Thus we have
\be 
	\mathbb{P}(j_t(X)|\mathscr{Y}_t) = U^*_t 
	\mathbb{Q}^t (X \vert \mathscr{C}_t ) U_t ,\qquad x\in\mathscr{B}
	\label{c-exp-P-Q}
\ee
where $\mathscr{C}_t =\mathrm{vN}(Z_s:0\le s\le t)$.  Note that 
$\mathscr{Y}_t=U_t^*\mathscr{C}_tU_t$ follows from the fact that
$U_s^*Z_sU_s=U_t^*Z_sU_t$ for $t\ge s$, the property we used in 
\S\ref{sec filtering measurements} to prove self-nondemolition of $Y_t$.
The ease with which we will now be able to manipulate $\mathbb{Q}^t(X 
\vert\mathscr{C}_t)$ highlights the usefulness of the transformation 
(\ref{c-exp-P-Q}).

Our strategy will be as follows.  We wish to calculate $\mathbb{Q}^t(X
\vert\mathscr{C}_t)$; however, the state $\mathbb{P}$ has the nice 
property that $Z_{s\le t}$, which generates $\mathscr{C}_t$, is a
$\mathbb{P}$-Wiener process.  We want to use the Bayes formula,
Lemma \ref{thm KS}, in order to express $\mathbb{Q}^t(X\vert\mathscr{C}_t)$
in terms of $\mathbb{P}$-conditional expectations.  We run into a problem, 
however, as the ``change of measure'' operator $U_t$ that relates 
$\mathbb{P}$ with $\mathbb{Q}^t$ does not satisfy the requirement of Lemma 
\ref{thm KS} that\footnote{
	If this were the case then we could calculate
	$Y_t=U_t^*Z_tU_t=Z_tU_t^*U_t=Z_t$, i.e.\ the observations
	would carry no information about the system and the filtering
	problem would be trivial.
} $U_t\in\mathscr{C}_t'$.  To solve this problem, we will replace $U_t$ by 
a different operator $V_t$ which is affiliated to $\mathscr{C}_t'$, but 
which still defines the same state in the sense that 
$\mathbb{P}(U_t^*XU_t)=\mathbb{P}(V_t^*XV_t)$ for every $X$.  The 
following technique, to our knowledge, first appeared in \cite{Hol90}; it 
replaces Girsanov's theorem in the quantum context.

\begin{lemma}
\label{lem:vtdilation}
	Let $V_t$ be the solution of the QSDE
	\begin{equation}
	\label{eq Vt}
		dV_t=\Big\{
			L\,(dA_t^*+dA_t)-\frac{1}{2}L^*L\,dt-iH\,dt
		\Big\}V_t.
	\end{equation}
	Then $V_t$ is affiliated to $\mathscr{C}_t'$ and 
	$\mathbb{Q}^t(X)=\mathbb{P}(V_t^*XV_t)$ for all 
	$X\in\mathscr{B}\otimes\mathscr{W}$.
\end{lemma}

A precise proof of this statement is not very insightful, see e.g.\
\cite{BH05} or \cite{Bel92a}.  However, it is not difficult to see why the
statement should be true.  Let us assume for simplicity that the state
$\rho$ on $\mathscr{B}$ is pure; we can always obtain a mixed state later
by taking convex combinations. Then $\mathbb{P}(X)=
\langle\psi\ten\Phi,X\,\psi\ten\Phi\rangle$ for some vector $\psi\in\qh$
(and $\Phi\in\qF$ is the vacuum vector).  To show that
$\mathbb{P}(U_t^*XU_t)=\mathbb{P}(V_t^*XV_t)$, it is thus sufficient to
show that $U_t\,\psi\ten\Phi=V_t\,\psi\ten\Phi$.  Now recall from 
$\S\ref{sec Fock}$ that any stochastic integral with respect to $A_t$ 
vanishes when it acts on the vacuum vector; hence
\begin{equation*}
\begin{split}
  U_t\,\psi\ten\Phi &=
	\left[I+\int_0^tLU_s\,dA_s^*-
		\int_0^tL^*U_s\,dA_s-
		\int_0^t\left(
			\frac{1}{2}L^*L+iH
			\right)U_s\,ds\right]\psi\ten\Phi
	\\
	&=
	\left[I+\int_0^tLU_s\,dA_s^*
		-\int_0^t\left(
			\frac{1}{2}L^*L+iH
			\right)U_s\,ds\right]\psi\ten\Phi,
\end{split}
\end{equation*}
and similarly we obtain for $V_t$ acting on the vacuum
$$
  V_t\,\psi\ten\Phi =
	\left[I+\int_0^tLV_s\,dA_s^*
		-\int_0^t\left(
			\frac{1}{2}L^*L+iH
			\right)V_s\,ds\right]\psi\ten\Phi.
$$
But as these expressions are the same, they should have the same solution
$U_t\,\psi\ten\Phi = V_t\,\psi\ten\Phi$.  In principle, we could change 
the integrand of the $A_t$-integral arbitrarily without affecting how the 
QSDE acts on the vacuum; in Lemma \ref{lem:vtdilation} we exploit this 
fact to modify $U_t$ precisely so that it is in the commutant of 
$\mathscr{C}_t$; indeed, Eq.\ \eqref{eq Vt} is driven only by the
noise $Z_t=A_t + A_t^*$ and its coefficients are in 
$\mathscr{B}\subset\mathscr{C}'_t$.

We are now ready to apply the Bayes formula, Lemma \ref{thm KS}.
Together with Lemma \ref{lem:vtdilation} and Eq.\ (\ref{c-exp-P-Q}),
we immediately obtain the following result.

\begin{theorem}[Noncommutative Kallianpur-Striebel]  \label{thm NC KS}
	Define for any system operator $X \in \mathscr{B}$ the 
	unnormalized conditional expectation
	\be
		\sigma_t(X)=U_t^*\,\mathbb{P}(V_t^*XV_t|\mathscr{C}_t)\,U_t
		\in \mathscr{Y}_t .
		\label{sigma-def}
	\ee
	Then the conditional expectation {\rm \er{jt-X-cdl-exp}} is given 
	by
	\begin{equation}
		\pi_t(X)=\frac{\sigma_t(X)}{\sigma_t(I)},\qquad
		\forall\,X\in\mathscr{B}.
		\label{nc ks}
	\end{equation}
\end{theorem}

We now obtain an explicit expression for $\sigma_t(X)$.  

\begin{theorem}[Unnormalized quantum filtering equation]
\label{thm DMZ equation}
The unnormalized conditional expectation $\sigma_t(X)$ satisfies the 
following linear QSDE:
\begin{equation}\label{eq Zakai}
	  d\sigma_t(X) = 
	  \sigma_t(\mathcal{L}_{L,H}(X))\,dt + 
	  \sigma_t(L^*X+XL)\,dY_t.
\end{equation}
\end{theorem}

To obtain (\ref{eq Zakai}) we will need to take conditional 
expectations of quantum It\^o integrals.  Let us briefly show how to do 
this.  First, we claim that if $K_t$ is an adapted process with $K_s$ 
affiliated to $\mathscr{C}_s'$, then
$\mathbb{P}(K_s|\mathscr{C}_t)=\mathbb{P}(K_s|\mathscr{C}_s)$ for $s \leq
t$. This follows from the fact that $\mathscr{C}_t=
\mathscr{C}_s\otimes\mathscr{C}_{[s,t]}$ and that $K_s$ is independent
from $\mathscr{C}_{[s,t]}$ by adaptedness.  Second, conditional
expectations and integrals can be exchanged as follows:
$$
                        \mathbb{P}\left(\left.
                                \int_0^tK_s\,ds\right|
                                \mathscr{C}_t\right)=
                        \int_0^t\mathbb{P}(K_s|\mathscr{C}_s)\,ds,
		\qquad
                        \mathbb{P}\left(\left.
                                \int_0^tK_s\,dZ_s\right|
                                \mathscr{C}_t\right)=
                        \int_0^t\mathbb{P}(K_s|\mathscr{C}_s)\,dZ_s.
$$
These properties are immediate if $K_t$ is a simple process, and a proof 
of the general case is not difficult.  

\begin{proof}
Using the quantum It\^o rules we have
$$
  V^*_tXV_t = X + \int_0^tV_s^*\mathcal{L}_{L,H}(X)V_s\,ds + 
	\int_0^tV^*_s(L^*X+XL)V_s\,d(A_s+A_s^*).
$$
We next take conditional expectations of the terms in this expression;
we obtain
\begin{multline*}
   	\BB{P}(V_t^*XV_t|\mathscr{C}_t) =  \BB{P}(X) + 
  	\int_0^t\BB{P}(V_s^*\mathcal{L}_{L,H}(X)V_s|\mathscr{C}_s)\,ds 
  \\
	+ \int_0^t\BB{P}(V^*_s(L^*X+XL)V_s|\mathscr{C}_s)\,d(A_s+A_s^*).
\end{multline*}
Another application of the quantum It\^o rules now yields \er{eq Zakai}.
\qquad
\end{proof}

By applying the It\^o rules to the noncommutative Kallianpur-Striebel 
formula \er{nc ks}, we obtain an expression for the normalized conditional 
state
\begin{equation}\label{eq qfilter}
	d\pi_t(X) = \pi_t(\mathcal{L}_{L,H}(X))dt + 
	\Big(\pi_t(L^*X+XL)-\pi_t(L^*+L)\,\pi_t(X)\Big)
	\Big(dY_t- \pi_t(L^*+L)\,dt\Big). 
\end{equation}
This (normalized) {\em quantum filtering equation} is a quantum analog
of the classical Kushner-Stratonovich equation of nonlinear filtering.
Note that this is a classical stochastic differential equation by the 
spectral theorem: it is a recursive equation that is only driven by the 
(commutative) observations $Y_t$.  Hence it can be implemented on a 
classical (digital) signal processor, as depicted in Fig.\ \ref{fig:model}.

\begin{remark}
Eq.\ (\ref{eq qfilter}) is expressed in terms of the conditional state 
$\pi_t(X)$, where $X\in\mathscr{B}$.  Now recall from \S\ref{sec BG} that 
any state on a finite-dimensional Hilbert space can be expressed as 
${\rm Tr}[\rho X]$ for some density matrix $\rho$.  Similarly, if $\qh$
(and hence $\mathscr{B}$) is finite-dimensional, then we can always 
write $\pi_t(X)={\rm Tr}[\rho_tX]$ where $\rho_t$, the conditional 
density matrix, is a (random) density matrix that is a function of the 
observations up to time $t$.  From Eq.\ (\ref{eq qfilter}) we obtain
explicitly
$$
	d\rho_t=-i[H,\rho_t]\,dt+
	(L\rho_tL^*-\tfrac{1}{2}L^*L\rho_t-\tfrac{1}{2}\rho_tL^*L)\,dt
	+(L\rho_t+\rho_tL^*-{\rm Tr}[(L+L^*)\rho_t]\rho_t)\,dW_t
$$
where $dW_t=dY_t-{\rm Tr}[(L+L^*)\rho_t]\,dt$.  In \S\ref{sec mtg} we will 
see that $W_t$ is a Wiener process.  It is this representation that is 
usually found in the physics literature. \qed
\end{remark}

\subsection{Photon counting measurements}
\label{sec misc photon}

We now consider the photon counting setup given by Eqs.\ (\ref{X}) and 
(\ref{YL}).  We would like to follow the same procedure as for homodyne 
detection.  The following lemma, which replaces Lemma 
\ref{lem:vtdilation}, suggests how to proceed.  The proof is identical to 
that of Lemma \ref{lem:vtdilation}.

\begin{lemma}
\label{lem:countvt}
	Let $U_t'$ be the solution of the QSDE
	$$
		dU_t'=\Big\{L'\,dA^*_t-L^{\prime*}\,dA_t-\frac{1}{2}
			L^{\prime*}L'\,dt-iH'\,dt
		\Big\}U_t'
	$$
	and let $V_t'$ be the solution of
	$$
		dV_t'=\Big\{L'(d\Lambda_t+dA^*_t+dA_t+dt)
			-\frac{1}{2}L^{\prime*}L'\,dt-L'\,dt-iH'\,dt
		\Big\}V_t'.
	$$
	Then $V_t'$ is affiliated to $\mathrm{vN}(\Lambda_s+A^*_s+A_s+s:
	s\le t)'$ and $\mathbb{P}({U_t'}^*XU_t')=
	\mathbb{P}({V_t'}^*XV_t')$.
\end{lemma}

Define $Z_t=\Lambda_t+A_t^*+A_t+t$ and $\mathscr{C}_t=\mathrm{vN}(Z_s:0\le
s\le t)$.  Lemma \ref{lem:countvt} directly provides us with a
nondemolition change of measure, provided that we rotate our problem so
that $\mathscr{Y}_t={U_t'}^*\mathscr{C}_tU_t'$ using a suitable unitary
operator $U_t'$.  Then, defining $\sigma_t(X)={U_t'}^*\,\mathbb{P}({V_t'}^*
XV_t'|\mathscr{C}_t)\,U_t'$, the Kallianpur-Striebel formula holds for
$\sigma_t(X)$.

Define $R_t$ as the solution of the QSDE
$$
	dR_t=(dA_t-dA_t^*-\tfrac{1}{2}dt)\,R_t
$$
Recall Example \ref{ex qsc}; evidently $R_t$ is a Weyl operator, 
and in particular $\Lambda_t=R_t^*Z_tR_t$.  But recall that 
$Y_t=U_t^*\Lambda_tU_t=U_t^*R_t^*Z_tR_tU_t$; thus $U_t'=R_tU_t$ is our 
rotation of choice.  Using the quantum It\^o rules we obtain
$$
	dU_t'=\Big\{
		(L-1)\,dA_t^* -(L^* -1)\,dA_t-\frac{1}{2}
		(L^* L+I-2L+2iH)\,dt
	\Big\}U_t',
$$
which corresponds to the nondemolition change of measure
$$
	dV_t'=\Big\{
		(L-1)\,dZ_t-\frac{1}{2}(L^*L-I+2iH)\,dt
	\Big\}V_t'.
$$
For $X\in\mathscr{B}$, we obtain using the quantum It\^o rules
\begin{equation*}
  d{V_t}'^*XV_t' = {V_t'}^*\big(\mathcal{L}_{L,H}(X)\big)V_t'\,dt + 
	{V_t'}^*(L^*XL-X)V_t'\,(dZ_t-dt). 
\end{equation*}  
Finally we obtain using the definition of $\sigma_t$ and the quantum 
It\^o rules
$$
  d\sigma_t(X) = \sigma_t(\mathcal{L}_{L,H}(X))\,dt + 
  \big(\sigma_t(L^*XL)-\sigma_t(X)\big)\big(dY_t-dt\big). 
$$
which is the unnormalized quantum filtering equation for counting 
observations.

Using the Kallianpur-Striebel formula $\pi_t(X)=\sigma_t(X)\,/\,\sigma_t(I)$ 
we can now obtain an expression for the normalized conditional state
\begin{equation*}
  d\pi_t(X) = \pi_t(\mathcal{L}_{L,H}(X))\,dt + 
  \left(\frac{\pi_t(L^*XL)}{\pi_t(L^*L)}-\pi_t(X)\right)
	\big(dY_t-\pi_t(L^*L)\,dt\big), 
\end{equation*}
which is the normalized quantum filtering equation for photon counting.


\section{The innovations method}
\label{sec mtg}

In this section we rederive the filtering equation for homodyne detection,
Eq.\ (\ref{eq qfilter}), using martingale methods that are analogous to
the classical case \cite{Bel92b,BGM04}. We follow the classical treatment
as in \cite{FKK72}, \cite[chapter 18]{RE82}, \cite[chapter 7]{WH85}.  
Martingale methods have enjoyed wide and successful application in
classical stochastic theory.  The procedure is less straightforward than
the reference probability method, however, and some familiarity with
classical filtering theory would be helpful (see e.g.\ \cite{DM81} for an
excellent introduction).

Let $\xi_t,\beta_t,\lambda_t,\mu_t$ be adapted processes affiliated to 
$\mathscr{Y}_t'$, where
\be
	\xi_t = \xi_0 + \int_0^t\beta_s\,ds + m_t =
	\xi_0 + \int_0^t\beta_s\,ds 
	+ \int_0^t (\lambda_s\,dA_s + \mu_s\,dA_s^*).
	\label{q-smtg}
\ee
The measurement process $Y_t$ is given by \er{Y}, and in what follows we 
write $h_t=j_t(L+L^*)$ and $Z_t=A_t+A^*_t$.  Note that the conditional
expectation $\hat\xi_t = \mathbb{P}(\xi_t \vert \mathscr{Y}_t)$ is well 
defined, and similarly for the coefficients $\beta_t,\lambda_t$ and $\mu_t$.

The main filtering result for a process of the form \er{q-smtg} is 
the following.

\begin{theorem}[Noncommutative Fujisaki-Kallianpur-Kunita]   
\label{thm smtg filter}
Under the above assumptions, the filtered process $\hat\xi_t$ satisfies 
the QSDE
\be
	d \hat \xi_t = \hat \beta_t\,dt + (\hat \lambda_t + \widehat{\xi_t
		h_t} - \hat \xi_t \hat h_t )\,dW_t
	\label{q-smtg-filter}
\ee
where $\hat r_t \equiv \mathbb{P}(r_t\vert\mathscr{Y}_t)$ for any $r_t$ 
affiliated to $\mathscr{Y}_t'$, and $dW_t = dY_t - \hat h_t\,dt$
defines the $\mathscr{Y}_t$-Wiener process (with respect to 
$\mathbb{P}$) $W_t$, called the {\em innovations process}.
\end{theorem}

The filtering expression \er{q-smtg-filter} is formally identical to the
classical case \cite[Theorem 18.11]{RE82}, \cite[Proposition 3.2]{WH85}.
Before we prove Theorem \ref{thm smtg filter}, we will show how to obtain 
the quantum filtering equation \er{eq qfilter} using this result.

\begin{corollary} \label{corol belavkin filter}
The conditional state $\pi_t(X)$ is given by Eq.\ {\rm \er{eq qfilter}}.
\end{corollary}

\begin{proof}
We set $\lambda_t=-j_t([X, L^*])$, $\mu_t=j_t([X, L])$,
$\beta_t = j_t(\mathcal{L}_{L,H}(X))$, and $\xi_t = j_t(X)$. Then 
$\widehat{\xi_t h_t} = \pi_t(X(L+L^*))$, $\hat\xi_t \hat h_t = \pi_t (X) 
\pi_t(L+L^*)$, $\hat\lambda_t = -\pi_t([X,L^*] )$, and $\hat\beta_t =
\pi_t(\mathcal{L}_{L,H}(X))$.  Hence using Eq.\ \er{q-smtg-filter}, 
Eq.\ \er{eq qfilter} follows.  \qquad 
\end{proof}

\begin{proof} {\em (Theorem \ref{thm smtg filter})}. 
{\em Step 1.} We first show that the process
$$
	M_t = \hat\xi_t - \hat\xi_0 - \int_0^t \hat\beta_s\,ds
$$
is a $\mathscr{Y}_t$-martingale, i.e.\ $\mathbb{P}(M_t\vert\mathscr{Y}_s)
= M_s$ for all $s \leq t$.  This property is equivalent to
$\mathbb{P}((M_t-M_s)K)=0$ for all $K\in \mathscr{Y}_s$, or
equivalently
$$
	\mathbb{P}\left[
	\left(\hat\xi_t - \hat\xi_s - \int_s^t\hat\beta_r\,dr\right) 
		K\right] = 
	\mathbb{P}\left[
	\left(\xi_t - \xi_s - \int_s^t\beta_r\,dr\right) 
		K\right] =
	\mathbb{P}[(m_t-m_s)K] = 0
$$
for all $K\in \mathscr{Y}_s$, where we have used Def.\ \ref{de 
conditional expectation} in the first step.  But as 
$K\in\mathscr{Y}_s\subset\mathscr{B}\otimes\mathscr{W}_{s]}$
$$
	\mathbb{P}[(m_t-m_s)K]=
	\mathbb{P}\left[
		K\int_s^t(\lambda_r\,dA_r+\mu_r\,dA_r^*)
	\right]=
	\mathbb{P}\left[
		\int_s^t(K\lambda_r\,dA_r+K\mu_r\,dA_r^*)
	\right]=0
$$
where we have used that the vacuum expectation of quantum It\^o integrals 
vanishes.  Thus we have demonstrated that $M_t$ is a 
$\mathscr{Y}_t$-martingale.

{\em Step 2.} We now show that $W_t$ is a Wiener process under 
$\mathbb{P}$.  We begin by verifying that the innovations process
\begin{equation}
	W_t = Y_t- \int_0^t \hat h_s\, ds
	\label{dy}
\end{equation}
is a $\mathscr{Y}_t$-martingale.  We need to show that
$\mathbb{P}[(W_t-W_s)K]=0$ for any $s\leq t$ and $K \in \mathscr{Y}_s$.  
This is equivalent to 
$$
	\mathbb{P}\left[\left(Y_t-Y_s-\int_s^t\hat h_r\, dr\right)K\right] 
	= \mathbb{P}\left[\left(Y_t-Y_s -\int_s^t h_r\,dr\right)K\right] 
	=0
$$
for all $K \in \mathscr{Y}_s$, where the second expression follows from the
definition of the conditional expectation.  But from \er{Y} we obtain
$$
	\mathbb{P}\left[\left(Y_t-Y_s -\int_s^t h_r\,dr\right)K\right]=
	\mathbb{P}[(A_t-A_s)K]=0
$$
as $K\in\mathscr{Y}_s\subset\mathscr{B}\otimes\mathscr{W}_{s]}$, 
$(A_t-A_s)\in\mathscr{W}_{[s,t]}$ and hence
$\mathbb{P}[(A_t-A_s)K]=\mathbb{P}(K)\,\mathbb{P}(A_t-A_s)=0$.
Thus $W_t$ is a $\mathscr{Y}_t$-martingale.

From \er{dy} we read off the It\^o rule $dW_t^2=dt$;  classically, a
process that obeys this property and is a martingale must be a Wiener
process by L{\'e}vy's Theorem (e.g. \cite[Lemma 18.7]{RE82}). But we can
simply apply the classical result, as $W_t$ is a commutative process (note
that $\hat h_t\in\mathscr{Y}_t$ for $s\le t$ by construction) and is hence
equivalent to the corresponding classical process obtained through the
spectral theorem.

Now that we have shown that $W_t$ is a Wiener process, we can try to 
represent the martingale $M_t$ as a stochastic integral with respect to 
$W_t$.  As usual in filtering theory the ordinary martingale 
representation theorem does not suffice for this purpose, but the 
representation theorem of Fujisaki-Kallianpur-Kunita (e.g.\ \cite[Theorem 
5.20]{LiS01}) allows us to conclude nonetheless that
\be
	M_t = \int_0^t \gamma_s\,dW_s\quad\Longrightarrow\quad
	\hat\xi_t=
	\hat\xi_0+\int_0^t \hat\beta_s\,ds+\int_0^t \gamma_s\,dW_s
	\label{M-dy}
\ee
for some adapted process $\gamma_t \in \mathscr{Y}_t$.

{\em Step 3.} 
We next obtain a first expression for $\widehat {\xi_tY_t}$:
\be
	\widehat {\xi_tY_t} = 
	\int_0^t [ \widehat{\beta_s Y_s} +
	\widehat{\xi_sh_s} + \hat\lambda_s ] ds + M_1(t) ,
	\label{term-1}
\ee
where $M_1(t)$ is a $\mathscr{Y}_t$-martingale.  As before, it suffices to 
show that 
$$
	\mathbb{P}((M_1(t)-M_1(s))K)=
	\mathbb{P}\left[\left({\xi_tY_t} -
	\int_0^t [ {\beta_s Y_s} +
	{\xi_sh_s} + \lambda_s ] ds\right)K
	\right]=0
$$
for all $K\in\mathscr{Y}_s$, where we have used the definition of the 
conditional expectation.  But
\begin{eqnarray*}
	d(\xi_t Y_t) &=& (d\xi_t) Y_t + \xi_t dY_t + d\xi_t dY_t
	\\
	&=&
	(\beta_t dt + dm_t) Y_t + \xi_t(hdt+dZ_t) + dm_t dZ_t
	\\
	&=&
	(\beta_t Y_t + \xi_t h_t + \lambda_t) dt + (Y_t \lambda_t + \xi_t) dA_t 
	+ (Y_t \mu_t +\xi_t) dA^\ast_t  .
\end{eqnarray*}
Hence exactly as before, it follows that $M_1(t)$ is a 
$\mathscr{Y}_t$-martingale.

{\em Step 4.} 
Next, we derive a second expression for $\widehat {\xi_tY_t}$:
\be
	\widehat {\xi_tY_t} = \int_0^t [\hat\beta_s Y_s + 
	\hat \xi_s \hat h_s + \gamma_s] ds + M_2(t) ,
	\label{term-2}
\ee
where $M_2(t)$ is a $\mathscr{Y}_t$-martingale.  To show this, note that 
$\widehat{\xi_tY_t} = \hat\xi_t Y_t$. By It\^o's rules,
\begin{eqnarray*}
	d(\hat\xi_t Y_t) &=& (d\hat\xi_t) Y_t + \hat\xi_t dY_t + d\hat\xi_t dY_t
	\\
	&=&
	(\hat\beta_t dt + \gamma_td W_t) Y_t + \hat\xi_t(\hat h_tdt+d W_t) 
		+ \gamma_t d W_t d W_t
	\\
	&=&
	(\hat\beta_t Y_t + \hat\xi_t \hat h_t + \gamma_t) dt + 
		(\gamma_t Y_t+\hat \xi_t) d W_t
\end{eqnarray*}
which establishes \er{term-2}.

{\em Step 5.}  
We can now identify  $\gamma_t$.  From \er{term-1} and \er{term-2} we have 
two representations for $\widehat {\xi_tY_t}$. By uniqueness, it follows 
that the finite variation terms are equal, viz.
$$
	\widehat{\beta_s Y_s} +
	\widehat{\xi_sh_s} + \hat\lambda_s =\hat\beta_s Y_s + \hat \xi_s 
	\hat h_s + \gamma_s .
$$
Therefore $\gamma_s = \widehat{\xi_sh_s} + \hat\lambda_s -  \hat \xi_s \hat 
h_s$ as required.
\qquad
\end{proof}

\bibliographystyle{siam}
\bibliography{ref}

\end{document}